\documentclass{article}
\usepackage{latexsym}
\usepackage{amssymb}
\usepackage{amsfonts}
\usepackage[ansinew]{inputenc}

\newtheorem{theorem}{Theorem}[section]
\newtheorem{lemma}{Lemma}[section]
\newtheorem{corollary}{Corollary}[section]
\newtheorem{proposition}{Proposition}[section]

\newtheorem{remark}{Remark}[section]

\newcommand{\pf}{{\bf Proof.~}} 

\def\Sp
 {\mathop{\rm Sp}\nolimits}

\def\Supp
 {\mathop{\rm Supp\,}\nolimits}

\def\lin
 {\mathop{\rm lin}\nolimits}

\def\Arg
 {\mathop{\rm Arg\,}\nolimits}

\def\rk
 {\mathop{\rm rank\,}\nolimits}

\def\Hom
 {\mathop{\rm Hom\,}\nolimits}

\def\codim
 {\mathop{\rm codim\,}\nolimits}

\def\re
 {\mathop{\rm Re\,}\nolimits}

\def\im
 {\mathop{\rm Im\,}\nolimits}

\def\Ch
 {\mathop{\rm Ch\,}\nolimits}

\def\Ex
 {\mathop{\rm Exp}\nolimits}
 
\def\card
 {\mathop{\rm card}\nolimits}

\begin{document}
\title{\bf Bochner transforms, perturbations and amoebae of holomorphic almost periodic mappings in tube domains}
\author{\scshape{Adelina Fabiano\thanks{Departement of Mathematics, University of Calabria, Arcavacata di Rende, Italy. {\tt fabiano@unical.it}}, Jacques Guenot\thanks{Departement of Mathematics, University of Calabria, Arcavacata di Rende, Italy. {\tt j.guenot@majise.it}}, James Silipo\thanks{Laboratory of Applications of Mathematics to Engineering, University of Calabria, Arcavacata di Rende, Italy. {\tt silipo@math.u-bordeaux1.fr}}}}
\date{}
%\markboth{J. Silipo et al.}{Amoebae of holomorphic almost periodic mappings}
\maketitle
\begin{abstract}
We give an alternative representation of the closure of the Bochner transform of a holomorphic almost periodic mapping in a tube domain. For such mappings we introduce a new notion of amoeba and we show that, for mappings which are regular in the sense of Ronkin, this new notion agrees with Favorov's one. We prove that the amoeba complement of a regular holomorphic almost periodic mapping, defined on~$\mathbb C^n$ and taking its values in~$\mathbb C^{m+1}$, is a Henriques~$m$-convex subset of~$\mathbb R^n$. Finally, we compare some different notions of regularity.
\end{abstract}

\section{Introduction and statement of the main result}

Let~$\Omega\subset\mathbb R^n$ be a non empty open convex set and let~$T_\Omega= \mathbb R^n+i\Omega\subset\mathbb C^n$ be the horizontal tube on the base~$\Omega$. A holomorphic mapping
$$
F=(f_1,\ldots,f_m):T_\Omega\longrightarrow\mathbb C^m
$$
is said to be {\sl almost periodic} in~$T_\Omega$ if for every~$\varepsilon>0$ and every~$\Omega^\prime\Subset\Omega$ there exists an~$n$-tuple~$t=t(\varepsilon,\Omega^\prime)$ of positive real numbers such that, for each~$\xi\in\mathbb R^n$, the~$n$-parallelotope~$\{x\in\mathbb R^n\mid \vert x_\ell-\xi_\ell\vert<t_\ell,\;1\leq\ell\leq n\}$ contains a point~$\tau$ for which
$$
\sup_{z\in T_{\Omega^\prime}}
\max_{1\leq\ell\leq m}
\vert f_\ell(z+\tau)-f_\ell(z)\vert<\varepsilon\,.
$$
Let~$HAP(T_\Omega,\mathbb C^m)$ denote the space of holomorphic almost periodic mappings defined on~$T_\Omega$ with values in~$\mathbb C^m$ endowed with the metrizable to\-po\-lo\-gy~$\tau(T_\Omega)$ of uniform convergence on every subtube~$T_{\Omega^\prime}$ with a relatively compact base~$\Omega^\prime\Subset\Omega$. 

Recall that by Bonchner's criterion, a holomorphic mapping
$$
F=(f_1,\ldots,f_m):T_\Omega\rightarrow\mathbb C^m
$$
is almost periodic if and only if the set
$$
{\mathcal B}(F)
=
\{(f_1(z+t),\ldots,f_m(z+t))\in HAP(T_\Omega,\mathbb C^m)\mid t\in\mathbb R^n\}\,,
$$
often referred to as the {\sl Bochner transform\/} of~$F$ (cf.~\cite{ap71}), is relatively compact. 
\par
\medskip
The {\sl amoeba}~${\mathcal F}_F$ (cf.~\cite{fa01}) of a mapping~$F\in HAP(T_\Omega,\mathbb C^m)$ is the closure of the image of its zero set~$V(F)$ under the projection~$\im:\mathbb C^n\rightarrow\mathbb R^n$ that maps each~$n$-tuple of complex numbers to the~$n$-tuple of their respective imaginary parts, i.e.
$$
{\mathcal F}_F=\overline{\im V(F)}\,.
$$
The notion of amoeba was originally introduced in 1994 by Gelfand, Kapranov and Zelevinski~\cite{gkz94} for algebraic subvarieties of tori. If~$V\subseteq(\mathbb C^*)^n$ is such a variety, the {\sl (polynomial) amoeba}~${\mathcal A}_V$ of~$V$ is the set
$$
{\mathcal A}_V
=
\{
(\log\vert z_1\vert,\ldots,\log\vert z_n\vert)\in\mathbb R^n
\mid
z\in V
\}\,.
$$
After their introduction, polynomial amoebae have been the object of intensive study. In 2001 this notion was adapted by Favorov~\cite{fa01} to the almost periodic setting and though this almost periodic version seems more natural than the usual one, the two notions have not yet been payed the same attention. 
\par
This paper is concerned with the study of Henriques' convexity of Favorov's amoeba complements. For the convenience of the reader we recall the notion of convexity introduced by Henriques~\cite{hen04}.
\par\medskip
Let~$m\in\mathbb N$,~$S\subseteq\mathbb R^n$ be an~$(m+1)$-dimensional oriented affine subspace and let~$Y\subseteq S$ be a subset. A class in the reduced homology group~${\widetilde H}_m(Y,\mathbb Z)$ is {\sl non negative} if its image in~${\widetilde H}_m(S\setminus\{x\},\mathbb Z)\simeq\mathbb Z$ is non negative, for every point~$x\in S\setminus Y$.
\par
A subset~$X\subseteq \mathbb R^n$ is~$m${\sl-convex} if, for every~$(m+1)$-dimensional oriented affine subspace~$S\subseteq\mathbb R^n$, the zero class is the only non negative class of the group~${\widetilde H}_m(S\cap X,\mathbb Z)$ that is sent to zero in~${\widetilde H}_m(X,\mathbb Z)$.
\par\medskip
In the paper~\cite{hen04}, Henriques proves the following theorem.
\begin{theorem}{\rm (Henriques \cite{hen04})}\label{henr}
Let~$V\subseteq(\mathbb C^*)^n$ be an~$(m+1)$-codimensional algebraic subvariety. Then, the amoeba complement~$\mathbb R^n\setminus{\mathcal A}_V$ is~$m$-convex.
\end{theorem}
In order to prove the almost periodic analogue of this result we need the following assumption of regularity introduced by Ronkin~\cite{ro90}.
\par
A mapping~$F\in HAP(T_\Omega,\mathbb C^m)$, with~$1\leq m\leq n$, is said to be {\sl regular} if, for every~$G\in\overline{{{\mathcal B}(F)}}$, the zero set of~$G$ is a complete intersection, (i.e. the zero set~$V(G)$ is either empty or~$m$-codimensional).
\par
We obtain the following main result, (cf. Corollary~\ref{buonis}).
\begin{theorem}
Let~$\Omega\subseteq\mathbb R^n$ be a non-empty, open, convex subset and let also~$F=(f_1,\ldots,f_{m+1})\in HAP(T_\Omega,\mathbb C^{m+1})$, with~$m+1\leq n$. If~$F$ is regular, then~$\Omega\setminus{\mathcal F}_F$ is an~$m$-convex subset.
\end{theorem}
This theorem shades some light on the topology of amoeba complements in the almost periodic setting and, in fact, it
generalizes a result shown in the paper~\cite{sil06} by the third author, where he considered the particular case of generic systems of exponential sums. As in that paper, the proof of the main result uses an alternative representation of the amoeba requiring a perturbation technique. 
\par
The sequel of the paper has the following structure. Section~\ref{prelim} fixes some notation and gives some preliminary information about holomorphic almost periodic mappings and Henriques' convexity, Section~\ref{perturbation} deals with a perturbation technique to be performed on holomorphic almost periodic mappings, Sections~\ref{Amoebae} and~\ref{Amoebacompl} are respectively devoted to amoebae and their complementary sets. Finally, in Section~\ref{fin} we compare some different conditions of regularity.

\section{Preliminaries}\label{prelim}
\subsection{Holomorphic almost periodic mappings in tube domains}

If~$F\in HAP(T_\Omega,\mathbb C^m)$, each component~$f_\ell$,~$1\leq\ell\leq m$, of~$F$ has a well-defined {\sl Bohr transform}
$$
a(\cdot,f_\ell):\mathbb R^n\longrightarrow\mathbb C\,,
$$ 
given, for every~$\lambda\in\mathbb R^n$, by
$$
a(\lambda,f_\ell)
=
\lim_{s\to+\infty}
{1\over (2s)^n}
\int_{\vert x_j\vert<s\,, j=1,\ldots,n}
e^{-i\langle x+iy,\lambda\rangle}
f_\ell(x+iy)\,dx\,.
$$
The Bohr transform of~$f_\ell$ is zero on the whole~$\mathbb R^n$ with, at most, the exception of a countable set called the {\sl spectrum} of~$f_\ell$,
$$
\Sp f_\ell
=
\{\lambda\in\mathbb R^n
\vert
a(\lambda,f_\ell)
\neq
0
\}\,.
$$
If~$G=(g_1,\ldots,g_m)\in\overline{{\mathcal B}(F)}$, then~$\Sp f_\ell=\Sp g_\ell$, for every~$1\leq\ell\leq m$.

By uniqueness theorem, if the Bohr transforms of two holomorphic almost periodic functions coincide, then the two functions themselves coincide.

Given~$F\in HAP (\mathbb C^n,\mathbb C^m)$, the~$\mathbb Z$-module~$\Xi_F$ of~$F$ (resp. the $\mathbb Q$-module $\lin_{\mathbb Q} \Xi_F$ of~$F$) is the additive subgroup (resp. the~$\mathbb Q$-linear subspace) of~$\mathbb R^n$ generated by the union of the spectra of the components of~$F$.

By an exponential sum we mean a finite non trivial~$\mathbb C$-linear combination of exponentials of the form~$e^{i\langle z,\lambda\rangle}$, with~$\lambda\in\mathbb R^n$. If~$f$ is such a function,~$f$ is holomorphic almost periodic and it has the following expression
$$
f(z)
=
\sum_{\lambda\in\Sp f}
a(\lambda,f)\,e^{i\langle z,\lambda\rangle}\,.
$$
If~$\Ex(T_\Omega,\mathbb C^m)$ denotes the subset of~$HAP(T_\Omega,\mathbb C^m)$ consisting of those mappings whose components are exponential sums, recall that, by Bochner-Fejer theorem,~$\Ex(T_\Omega,\mathbb C^m)$ is dense in~$HAP(T_\Omega,\mathbb C^m)$.
\par
Given~$F=(f_1,\ldots,f_m)\in HAP(T_\Omega,\mathbb C^m)$, though there may exist several different sequences of mappings belonging to~$\Ex(T_\Omega,\mathbb C^m)$ and converging to~$F$, some of them, introduced by Bochner, have a particularly nice form. 
For fixed~$1\leq\ell\leq m$, let~$\{\omega_r\}_{r\in\mathbb N}$ be a fixed basis of~$\lin_{\mathbb Q}\Xi_F$. For any~$j\in\mathbb N$, consider the exponential sum
$$
Q_j(f_\ell)
=
\sum
\prod_{r=1}^j
\bigg(
1-{\vert \nu_r\vert\over(j!)^2}
\bigg)
a\bigg(
\sum_{r=1}^j \nu_r(\omega_r/j!),f_\ell
\bigg)
\exp
\bigg({i\sum_{r=1}^j\nu_r\langle z, \omega_r/j!\rangle}\bigg)
$$
where the first sum is taken on the set of~$j$-tuples of integers~$\nu_1,\ldots,\nu_j$ having absolute values smaller or equal to~$(j!)^2$. Of course, only the~$j$-tuples such that
$$
\sum_{s=1}^j \nu_s{\omega_s\over j!}\in\Sp f_\ell
$$
contribute to the sum and, if~$\lambda=\sum_{s=1}^j \nu_s(\omega_s/ j!)\in\Sp f_\ell$, the corresponding product
$$
\mu_{\lambda,j}
=\prod_{r=1}^j
\bigg(
1-{\vert \nu_r\vert\over(j!)^2}
\bigg)
$$
is called the~$j$-th multiplier of~$f_\ell$ relative to~$\lambda$. For any~$\lambda\in\Sp f_\ell$, one has
$$
0\leq\mu_{\lambda,j}\leq 1\,,
$$
for every~$j\in\mathbb N$ and also
\begin{eqnarray}\label{multipl}
\lim_{j\to+\infty}
\mu_{\lambda,j}
=1\,.
\end{eqnarray}

In the sequel, for any holomorphic almost periodic mapping~$F$, with a fixed basis~$\{\omega_r\}_{r\in\mathbb N}$ of~$\lin_{\mathbb Q}\Xi_F$, the sequence of mappings
$$
Q_j(F)=((Q_j(f_1),\ldots,Q_j(f_m))\,,\;j\in\mathbb N\,,
$$
will be referred to as {\sl the Bochner-Fejer approximation of}~$F$, corresponding to the chosen basis~$\{\omega_r\}_{r\in\mathbb N}$. The choice of the basis does not appear in the notation~$Q_j(F)$ because we will never consider Bochner-Fejer approximations corresponding to different choices of basis.
\par

\subsection{Henriques' convexity}

If~$m\geq n-1$, every subset of~$\mathbb R^n$ is~$m$-convex, therefore, when dealing with~$m$-convexity in~$\mathbb R^n$, we will always assume~$m<n-1$. In this case, we may regard the~$m$-convexity of a subset~$X\subseteq\mathbb R^n$ as a statement of partial injectiveness for the morphisms~$\iota_{_{S,X}}:{\widetilde H}_m(S\cap X)\rightarrow{\widetilde H}_m(X)$ induced in reduced homology by the inclusions~$S\cap X\subseteq X$, as~$S$ varies among the oriented affine~$(m+1)$-dimensional subspaces of~$\mathbb R^n$. In fact, if~$X$ is~$m$-convex, the restriction of such a morphism to the submonoid~${\widetilde H}_m^+(S\cap X)$ of non negative classes is really injective. For~$m=0$ this partial injectiveness of the morphism~$\iota_{_{S,X}}$ is equivalent to its full injectiveness and, though we cannot present any counterexample yet, we expect that this is no longer the case for~$m>0$. However, we report that, when~$X$ is the amoeba complement of an~$(m+1)$-codimensional algebraic subvariety of~$(\mathbb C^*)^n$, Henriques~\cite{hen04} and Mikhalkin~\cite{mi04} expect the morphisms~$\iota_{_{S,X}}$ to be actually injective. As far as we know, no proof for this statement is presently available.
\par
We still recall that an open subset of~$\mathbb R^n$ is~$0$-convex if and only if each of its connected components is convex. Further, an~$m$-convex subset is also~$(m+1)$-convex (cf.~\cite{sil06}, Lemma~4.1), so that the weakness of this condition increases with~$m$ until it becomes the empty condition when~$m\geq(n-1)$. Observe finally that~$m$-convexity is preserved by linear authomorphismes of~$\mathbb R^n$ with positive determinants, (cf.~\cite{sil06}, Lemma~4.2).
\par
 
\subsection{Some useful lemmas}
Before embarking on the technical core of the paper, we collect here some easy or well-known lemmas which will be frequently used in the rest of the paper.
\begin{lemma}\label{krom}
Let~$\{\omega_\ell\}_{\ell\in\mathbb N}\subset\mathbb R^n$ a sequence of~$\mathbb Q$-linearly independent elements. Then the additive subgroup
$$
\mathbb G
=
\{x\in\mathbb R^{\mathbb N}\mid
x_\ell=\langle t,\omega_\ell\rangle+2\pi p_\ell\,,
\hbox{\rm with}\;
t\in\mathbb R^n\;
\hbox{\rm and}\;
p_\ell\in\mathbb Z\;
\hbox{\rm for every}\;\ell\in\mathbb N\}
$$
is dense in the product space~$\mathbb R^{\mathbb N}$.
\end{lemma}
\pf
Let~$r\in\mathbb N$ and let us consider the additive subgroup
$$
\mathbb G_r
=
\{x\in\mathbb R^r\mid
x_\ell=\langle t,\omega_\ell\rangle+2\pi p_\ell\;,
\hbox{\rm with}\;
t\in\mathbb R^n\;
\hbox{\rm and}\;
p_\ell\in\mathbb Z\;
\hbox{\rm for}\;1\leq\ell\leq r\}.
$$
The product topology on~$\mathbb R^{\mathbb N}$ induces the ordinary topology on~$\mathbb R^r$ and, by a Kronecker type result, (a proof of which is given for example in~\cite{sil06}, Theorem~3.1),~$\mathbb G_r$ is dense in~$\mathbb R^r$ with respect to this topology, so
$$
\overline\mathbb G
=
\overline{\bigcup_{r\in\mathbb N}\mathbb G_r}
\supseteq
\bigcup_{r\in\mathbb N}\overline\mathbb G_r
=
\bigcup_{r\in\mathbb N}\mathbb R^r
=
\mathbb R^{\mathbb N}\,,
$$
which proves the lemma.\hfill$\square$

\par

\begin{lemma}\label{comp}
Let~$\mathbb G$ any abelian group, then its Pontryagin dual group
$$
\Ch\mathbb G=\Hom_{\mathbb Z}(\mathbb G,\mathbb S^1)
$$
is compact in the weak topology.
\end{lemma}
\pf
The space~$(\mathbb S^1)^{\mathbb G}$ of all mappings from~$\mathbb G$ to~$\mathbb S^1$ is compact in the weak (i.e. product) topology.~$\Ch\mathbb G$ is a closed subspace, in fact, if~$(\chi_s)\subset\Ch\mathbb G$ is a sequence converging to~$\chi\in(\mathbb S^1)^{\mathbb G}$, then
$$
\chi(\lambda\mu)
=
\lim_{s\to+\infty} \chi_s(\lambda\mu)
=
\lim_{s\to+\infty} \chi_s(\lambda)\chi_s(\mu)
=
\chi(\lambda)\chi(\mu)\,,
$$
for every~$\lambda,\mu\in\mathbb G$.\hfill$\square$

\par
Finally, we recall the following result, whose proof is a special case of that of Lemma~3.2 in Rashkovskii's paper~\cite{ras01}, to which the reader is referred.
\begin{lemma}{\rm (Rashkovskii \cite{ras01})}\label{ras}
Let~$F\in HAP(T_\Omega,\mathbb C^m)$ be regular and let also $(G_j)\subset HAP(T_\Omega,\mathbb C^m)$ be a sequence converging to~$F$ in the topology~$\tau(T_\Omega)$. Then, for every~$\Omega^\prime\Subset\Omega$ there exists a~$j_0\in\mathbb N$ such that the mapping~$G_j$ is regular on~$T_{\Omega^\prime}$ as soon as~$j\geq j_0$.
\end{lemma}
\par

\section{Perturbation by characters}\label{perturbation}

In this section we consider an alternative description of the closure of the Bochner transform of a holomorphic almost periodic mapping. When dealing with mappings having exponential sums as components, the description given in Proposition~\ref{facilis} below is more effective than the usual one. The key ingredient is a perturbation technique that goes back to some early works of Alain Yger on spectral synthesis for convolution equations (cf.~\cite{y79},~\cite{yg79}).

If~$F=(f_1,\ldots,f_m)\in\Ex(T_\Omega,\mathbb C^m)$ and~$\chi\in\Ch\mathbb R^n$, define the {\sl perturbed mapping}
$$
F_\chi=((f_1)_\chi,\ldots,(f_m)_\chi)
$$ 
as the mapping with the~$\ell$-th component~$(f_\ell)_\chi$,~$1\leq \ell\leq m$, given by
$$
(f_\ell)_\chi(z)
=
\sum_{\lambda\in\Sp f_\ell}
a(\lambda,f_\ell)\,
\chi(\lambda)\,
e^{i\langle z,\lambda\rangle}\,.
$$
The perturbing mapping
$$
\Ch\mathbb R^n\times\Ex(T_\Omega,\mathbb C^m)\longrightarrow\Ex(T_\Omega,\mathbb C^m)\,,
$$
that associates to a given character~$\chi$ and a given mapping~$F$ the perturbed mapping~$F_\chi$, is clearly an action.
\par
\begin{proposition}\label{facilis}
The perturbing mapping is continuous with respect to the first argument.
\end{proposition}
\pf
Fix~$F=(f_1,\ldots,f_m)\in\Ex(T_\Omega,\mathbb C^m)$ and let~$(\chi_j)\subset\Ch\mathbb R^n$ be a sequence converging to a limiting character~$\chi\in\Ch\mathbb R^n$.
Let, for every fixed~$\Omega^\prime\Subset\Omega$,
$$
M_{\Omega^\prime}
=
\max_{1\leq\ell\leq m}
\bigg[
\bigg(
\max_{\lambda\in\Sp f_\ell}
\vert 
a(\lambda,f_\ell)
\vert
\bigg)
\bigg(
\max_{z\in T_{\Omega^\prime}}
e^{-\im\langle z,\lambda\rangle}
\bigg)
\bigg]\,,
$$
then
$$
\sup_{z\in T_{\Omega^\prime}}
\vert
(f_\ell)_\chi(z)-(f_\ell)_{\chi_j}(z)
\vert
\leq
M_{\Omega^\prime}
\sum\limits_{\lambda\in\Sp f_\ell}
\vert
\chi(\lambda)-\chi_j(\lambda)
\vert\,.
$$
It follows that~$F_\chi=\lim_{j\to+\infty} F_{\chi_j}$ in the topology~$\tau(T_\Omega)$.\hfill$\square$
\par
\begin{remark}
{\rm
Observe that, though the perturbing mapping described above is continuous with respect to its first argument, it is not continuous with respect to its second argument. The reason is that the characters are not supposed to be continuous in the ordinary topology of~$\mathbb R^n$. For example, let~$\gamma\in\mathbb R\setminus\mathbb Q$, let~$(\lambda_j)\subset\mathbb Q$ be a sequence converging to~$\gamma$ and let~$\chi$ be a character of~$\mathbb R$ such that~$\chi(\lambda)=1$, for every~$\lambda\in\mathbb Q$, and~$\chi(\gamma)=i$. Then, the sequence~$(e^{i\lambda_j z})\subset HAP(\mathbb C,\mathbb C)$ coincides with its perturbation by~$\chi$, but this is not the case for its limit~$e^{i\gamma z}$, whose perturbation is equal to~$ie^{i\gamma z}$. Consequently, the perturbed sequence does not converge to the perturbed limit and this is not a surprise since the character~$\chi$ is not continuous in the ordinary topology of~$\mathbb R$.}
\end{remark}
\par
\begin{lemma}\label{figo}
If~$F=(f_1,\ldots,f_m)\in \Ex(T_\Omega,\mathbb C^m)$, then the orbit
$$
\{F_\chi\in\Ex(T_\Omega,\mathbb C^m)\mid\chi\in\Ch\mathbb R^n\}
$$
of~$F$ equals the closure~$\overline{{\mathcal B}(F)}$, in~$HAP(T_\Omega,\mathbb C^m)$ with the topology~$\tau(T_\Omega)$, of its Bochner transform.
\end{lemma}
\pf
If~$F$ is a constant mapping the theorem is trivial, so let~$F$ be non constant.
We firstly show that~$\overline{{\mathcal B}(F)}$ is included in the orbit.
\par 
Let~$(t_s)\subset\mathbb R^n$ be a sequence such that~$F(z+t_s)$ converges in the to\-po\-lo\-gy~$\tau(T_\Omega)$ and let~$G=(g_1,\ldots,g_m)$ be the limiting mapping. For~$1\leq \ell\leq m$, we have
$$
f_\ell(z)
=
\sum_{\lambda\in\Sp f_\ell}
a(\lambda,f_\ell)\,
e^{i\langle z,\lambda\rangle}\,,
$$
so, for every~$s\in\mathbb N$, we have
$$
f_\ell(z+t_s)
=
\sum_{\lambda\in\Sp f_\ell}
a(\lambda,f_\ell)\,
e^{i\langle t_s,\lambda\rangle}
e^{i\langle z,\lambda\rangle}
=
(f_\ell)_{\chi_s}(z)\,,
$$
where~$\chi_s\in\Ch\mathbb R^n$ denotes the character given by~$\chi_s(\cdot)=e^{i\langle t_s,\cdot\rangle}$. By compactness~$(\chi_s)$ has a convergent subsequence~$(\chi_{s_q})$; if~$\chi$ is its limit, then, by~Poposition~\ref{facilis},~$G=\lim_{q\to+\infty} F_{\chi_{s_q}}=F_\chi$ in the topology~$\tau(T_\Omega)$.
\par
We now prove the other inclusion.
Let~$\chi\in\Ch\mathbb R^n$ and~$\omega_1,\ldots,\omega_r$ be a system of free generators of~$\Xi_F$, then, for any~$1\leq k\leq m$, we have
$$
f_k(z)
=
\sum_{\mu\in\mathbb Z^r}
a_{k,\mu}
\prod_{\ell=1}^r
e^{i\mu_\ell\langle z,\omega_\ell\rangle}\,,
$$
where all but a finite number of~$a_{k,\mu}$'s are different from zero. Set, for any index~$1\leq\ell\leq r$,~$\theta_\ell=\Arg\chi(\omega_\ell)$, hence
$$
(f_k)_\chi(z)
=
\sum_{\mu\in\mathbb Z^r}
a_{k,\mu}
\prod_{\ell=1}^r
e^{i\mu_\ell(\theta_\ell+\langle z,\omega_\ell\rangle)}\,,
$$
for any~$1\leq k\leq m$, and Proposition~\ref{krom} implies the existence of a sequence~$(t_s)\subset\mathbb R^n$ such that
$$
\lim_{s\to+\infty}
\langle t_s,\omega_\ell\rangle=\theta_\ell
\qquad
\hbox{\rm mod}\;2\pi\mathbb Z\,,
$$
hence
$$
\lim_{s\to+\infty}f_k(z+t_s)=(f_k)_\chi(z)
$$ 
uniformly on every compactly based subtube of~$T_\Omega$, for every~$1\leq k\leq m$, that is
$$
\lim_{s\to+\infty}F(z+t_s)=F_\chi(z)
$$
in the topology~$\tau(T_\Omega)$. \hfill$\square$
\par
\medskip
The next theorem shows how to pass to the general case.
\begin{theorem}\label{defpert}
Let~$F=(f_1,\ldots,f_m)\in HAP(T_\Omega,\mathbb C^m)$ and let~$\chi\in\Ch\mathbb R^n$. If~$\{\omega_r\}_{r\in\mathbb N}$ is a basis of~$\lin_\mathbb Q\Xi_F$ and~$Q_j(F)\subset\Ex(T_\Omega,\mathbb C^m)$ is the corresponding Bochner-Fejer approximation of~$F$, then the perturbed sequence~$(Q_j(F))_\chi$ admits a limit~$F_\chi$.
\end{theorem}
\pf
We show that the perturbed sequence is a Cauchy sequence in the complete space~$HAP(T_\Omega,\mathbb C^m)$ and we do it by contradiction. Assume there exits an~$\Omega^\prime\Subset\Omega$ and an~$\varepsilon>0$ such that for any~$\nu\in\mathbb N$ there are an integer~$j_\nu>\nu$, an integer~$p_\nu$ and a point~$z_\nu\in T_{\Omega^\prime}$ providing the relation
\begin{eqnarray}\label{a}
\max_{1\leq \ell\leq m}
\vert
(Q_{j_\nu+p_\nu}(f_\ell))_\chi(z_\nu)-(Q_{j_\nu}(f_\ell))_\chi(z_\nu)\vert
>\varepsilon\,.
\end{eqnarray}
We prove that this assumption leads to a contradiction. In fact let, for any~$r,j\in\mathbb N$,~$\theta_{r,j}=\Arg\chi(\omega_r/j!)$ and let, by~Lemma~\ref{krom},~$(t_{s,j})\subset\mathbb R^n$ be a sequence such that
$$
\lim_{s\to+\infty}
\langle (\omega_r/j!),t_{s,j}\rangle
=
\theta_{r,j}\,,
\qquad
{\rm mod}
\,2\pi\mathbb Z\,,
$$
for any~$r,j\in\mathbb N$. If~$r,j,p\in\mathbb N$ are fixed,
$$
\chi(\omega_r/j!)
=
\chi(\omega_r/(j+p)!)^{(j+p)!/j!}
\,,
$$
so~$\theta_{r,j}=\theta_{r,j+p}(j+p)!/j!$ and
$$
\lim_{s\to+\infty}
\langle (\omega_r/j!),t_{s,j}\rangle
=
\lim_{s\to+\infty}
\langle (\omega_r/j!),t_{s,j+p}\rangle
\qquad
{\rm mod}
\,2\pi\mathbb Z\,.
$$
The expression of~$Q_j(F)$ shows that the~$\mathbb Z$-module of~$Q_j(F)$ is freely generated by finitely many vectors of the form~$(\omega_r/j!)$,~$r\in\mathbb N$, so (as in the proof of Lemma~\ref{figo}), for every~$j,p\in\mathbb N$,
\begin{equation}\label{compatib}
(Q_j(F))_\chi(z)
=
\lim_{s\to+\infty}
Q_j(F)(z+t_{s,j})
=
\lim_{s\to+\infty}
Q_j(F)(z+t_{s,j+p})\,,
\end{equation}
in the topology~$\tau(T_\Omega)$. Now, the sequence~$Q_j(F)$ is Cauchy, so there is an index~$\nu_1=\nu_1(\varepsilon,\Omega^\prime)\in\mathbb N$ such that, for any integer~$\nu>\nu_1$,
$$
\sup_{z\in T_{\Omega^\prime}}
\max_{1\leq \ell\leq m}
\vert
Q_{j_\nu+p_\nu}(f_\ell)(z)-Q_{j_\nu}(f_\ell)(z)\vert
<(\varepsilon/3)\,,
$$
in particular
\begin{equation}\label{b}
\max_{1\leq \ell\leq m}
\vert
Q_{j_\nu+p_\nu}(f_\ell)(z_\nu+t_{s,j_\nu+p_\nu})-Q_{j_\nu}(f_\ell)(z_\nu+t_{s,j_\nu+p_\nu})\vert
<(\varepsilon/3)\,,
\end{equation}
for any~$\nu>\nu_1$. Also, by the relations~(\ref{compatib}), when~$\nu>\nu_1$, there is an integer~$s_1$, (depending on~$\varepsilon$ and on~$\nu$), such that, for any integer~$s>s_1$, we have
\begin{equation}\label{c}
\max_{1\leq \ell\leq m}
\vert
(Q_{j_\nu+p_\nu}(f_\ell))_\chi(z_\nu)-(Q_{j_\nu+p_\nu}(f_\ell))(z_\nu+t_{s,j_\nu+p_\nu})\vert
<(\varepsilon/3)\,
\end{equation}
and
\begin{equation}\label{d}
\max_{1\leq \ell\leq m}
\vert
(Q_{j_\nu}(f_\ell))(z_\nu+t_{s,j_\nu+p_\nu})-(Q_{j_\nu}(f_\ell))_\chi(z_\nu)
\vert
<(\varepsilon/3)\,.
\end{equation}
In conclusion, if~$\nu>\nu_1$ and~$s>s_1$, the relations~(\ref{b}),~(\ref{c}) and~(\ref{d}) imply that
$$
\max_{1\leq \ell\leq m}
\vert
(Q_{j_\nu+p_\nu}(f_\ell))_\chi(z_\nu)-(Q_{j_\nu}(f_\ell))_\chi(z_\nu)\vert
<\varepsilon\,,
$$
which contradicts the relation~(\ref{a}).\hfill$\square$

\par
\medskip
For every~$F\in HAP(T_\Omega,\mathbb C^m)$ and any~$\chi\in\Ch\mathbb R^n$, according to Theorem~\ref{defpert}, we set by definition
$$
F_\chi
=
\lim_{j\to+\infty}
Q_j(F)_\chi\,,
$$
where the limit is understood in the topology~$\tau(T_\Omega)$ and~$Q_j(F)$ is the Bochner-Fejer approximation of~$F$ corresponding to some choice of a basis for~$\lin_{\mathbb Q}\Xi_F$.
\par
\medskip
From now on, given~$F\in HAP(T_\Omega,\mathbb C^m)$, we will denote by~${\mathcal O}(F)$ the orbit of the mapping~$F$ with respect to the action of~$\Ch\mathbb R^n$. 
\par
\begin{remark}\label{importante}
{\rm
Observe that the perturbation of a fixed~$F\in HAP(T_\Omega,\mathbb C^m)$ may be performed using the $\mathbb S^1$-characters of any other abelian group~$\mathbb G$ containing the~$\mathbb Z$-module of~$F$. In fact, as the group~$\mathbb S^1$ is divisible, it is an injective object in the category of abelian groups, (i.e. the restriction map~$\Ch\mathbb G\rightarrow\Ch\Xi_F$ is surjective), so
$$
{\mathcal O}(F)
=
\{F_\chi\mid \chi\in\Ch\Xi_F\}
=
\{F_\chi\mid \chi\in\Ch\mathbb G\}\,.
$$
}
\end{remark}
\par
\begin{theorem}\label{figo2}
If~$F\in HAP(T_\Omega,\mathbb C^m)$, then~${\mathcal O}(F)=\overline{{\mathcal B}(F)}$.
\end{theorem}
\pf
Fix a basis of~$\lin_{\mathbb Q}\Xi_F$ and let also~$\chi\in\Ch\mathbb R^n$ be fixed. The proof of~Theorem~\ref{defpert} implies the existence of a sequence~$(t_{s,j})\subset\mathbb R^n$ such that, for every~$\Omega^\prime\Subset\Omega$,
$$
F_\chi(z)
=
\lim_{j\to+\infty}
Q_j(F)_\chi(z)
=
\lim_{j\to+\infty}
\lim_{s\to+\infty}
Q_j(F)(z+t_{s,j})\,,
$$
uniformly on~$T_{\Omega^\prime}$. For any fixed~$s\in\mathbb N$, let~$(t_{s,j_q})_{q\in\mathbb N}$ be such that the sequence~$F(z+t_{s,j_q})$ converges in the topology~$\tau(T_\Omega)$ and let
$$
G^{(s)}=(g_1^{(s)},\ldots, g_m^{(s)})
$$
denote the limiting mapping. We affirm that, for any~$s\in\mathbb N$,
$$
\lim_{q\to+\infty}
Q_{j_q}(F)(z+t_{s,j_q})
=
G^{(s)}(z)
$$
in the topology~$\tau(T_\Omega)$. In fact, suppose by contradiction that there is an~$\tilde s\in\mathbb N$, an~$\Omega^\prime\Subset\Omega$ and an~$\varepsilon>0$ such that, for every~$\nu\in\mathbb N$ we may find a~$q>\nu$ and a~$z_q\in T_{\Omega^\prime}$ for which
\begin{equation}\label{star}
\max_{1\leq\ell\leq m}
\vert
g_\ell^{(\tilde s)}(z_q)-Q_{j_q}(f_\ell)(z_q+t_{\tilde s,j_q})\vert
>\varepsilon\,.
\end{equation}
We already know that there exists a~$\tilde\nu\in\mathbb N$, (depending on~$\tilde s,\Omega^\prime$ and~$\varepsilon$), such that
$$
\sup_{z\in T_{\Omega^\prime}}
\max_{1\leq\ell\leq m}
\vert
g_\ell^{(\tilde s)}(z)-f_\ell(z+t_{\tilde s,j_q})\vert
<(\varepsilon/2)
$$
and
$$
\sup_{z\in T_{\Omega^\prime}}
\max_{1\leq\ell\leq m}
\vert
f_\ell(z)-Q_{j_q}(f_\ell)(z)\vert
<(\varepsilon/2)
$$
as soon as~$q>\tilde\nu$, so, when~$q>\tilde\nu$, we have
$$
\max_{1\leq\ell\leq m}
\vert
g_\ell^{(\tilde s)}(z_q)-f_\ell(z_q+t_{\tilde s,j_q})\vert
<(\varepsilon/2)
$$
and
$$
\max_{1\leq\ell\leq m}
\vert
f_\ell(z_q+t_{\tilde s,j_q})-Q_{j_q}(f_\ell)(z_q+t_{\tilde s,j_q})\vert
<(\varepsilon/2)\,.
$$
Consequently
$$
\max_{1\leq\ell\leq m}
\vert
g_\ell^{(\tilde s)}(z_q)
-
Q_{j_q}(f_\ell)(z_q+t_{\tilde s,j_q})
\vert
<\varepsilon\,,
$$
which is in contradiction with~(\ref{star}). Finally, for any~$s\in\mathbb N$,~$G^{(s)}\in\overline{{\mathcal B}(F)}$ (by definition), so by compactness there exists a subsequence~$(s_\iota)_{\iota\in\mathbb N}$ such that~$G^{(s_\iota)}$ converges in the topology~$\tau(T_\Omega)$. Hence, by uniform convergence,
\begin{eqnarray*}
F_\chi(z)
&=&
\lim_{q\to+\infty}
\lim_{\iota\to+\infty}
Q_{j_q}(F)(z+t_{s_\iota,j_q})
\\
&=&
\lim_{\iota\to+\infty}
\lim_{q\to+\infty}
Q_{j_q}(F)(z+t_{s_\iota,j_q})
\\
&=&
\lim_{\iota\to+\infty}
G^{(s_\iota)}(z)\,,
\end{eqnarray*}
in the topology~$\tau(T_\Omega)$. By the arbitrary choice of~$\chi\in\Ch\mathbb R^n$, it follows that~${\mathcal O}(F)\subseteq\overline{{\mathcal B}(F)}$.
\par
Now let~$G\in\overline{{\mathcal B}(F)}$ and let~$(t_s)\subset\mathbb R^n$ be a sequence such that
$$
G(z)=\lim_{s\to+\infty}F(z+t_s)
$$
uniformly in~$T_{\Omega^\prime}$, for every~$\Omega^\prime\Subset\Omega$. If~$(F_j)\subset\Ex(T_\Omega,\mathbb C^m)$ is the Bochner-Fejer approximation of~$F$ corresponding to some choice of a basis for~$\lin_{\mathbb Q}\Xi_F$, then we have
$$
G(z)
=
\lim_{s\to+\infty}F(z+t_s)
=
\lim_{s\to+\infty}
\lim_{j\to+\infty}
F_j(z+t_s)
$$
uniformly in~$T_{\Omega^\prime}$, for every~$\Omega^\prime\Subset\Omega$. If, for every~$s\in\mathbb N$,~$\chi_s\in\Ch\mathbb R^n$ denotes the character given by~$\chi_s(\cdot)= e^{i\langle t_s,\cdot\rangle}$, then
we have
$$
F_j(z+t_s)
=
(F_j)_{\chi_s}(z)\,,
$$
for every~$j,s\in\mathbb N$. By compactness, the sequence~$(\chi_s)$ has a convergent subsequence~$(\chi_{s_q})$. Let~$\chi$ be the limit of this subsequence, then, by uniform convergence on every compactly based subtube of~$T_\Omega$,
$$
G(z)
=
\lim_{q\to+\infty}
\lim_{j\to+\infty}
(F_j)_{\chi_{s_q}}(z)
=
\lim_{j\to+\infty}
(F_j)_{\chi}(z)
=
F_\chi(z)
$$
in the topology~$\tau(T_\Omega)$, whence~$G\in{\mathcal O}(F)$. The theorem is thus proved.\hfill$\square$
\par
\medskip
\begin{remark}
{\rm 
Theorem~\ref{figo2} easily implies that a mapping~$F\in HAP(T_\Omega,\mathbb C^m)$, with~$1\leq m\leq n$, is regular if and only if~$V(F_\chi)$ is a complete intersection, for every $\chi\in\Ch\mathbb R^n$.}
\end{remark}
\begin{corollary}\label{utile1}
Let~$F=(f_1,\ldots,f_m)\in HAP(T_\Omega,\mathbb C^m)$ and let~$Q_j(F)$ be the Bocner-Fejer approximation of~$F$ corresponding to a fixed basis of~$\lin_{\mathbb Q}\Xi_F$. Then, for  every~$\Omega^\prime\Subset\Omega$, every~$j\in\mathbb N$ and every character~$\chi\in\mathbb R^n$,
\begin{eqnarray}\label{magnifica}
\sup_{z\in T_{\Omega^\prime}}
\max_{1\leq\ell\leq m}
\vert
Q_j(f_\ell)_\chi(z)-(f_\ell)_\chi(z)
\vert
=
\sup_{z\in T_{\Omega^\prime}}
\max_{1\leq\ell\leq m}
\vert
Q_j(f_\ell)(z)-f_\ell(z)
\vert\,.
\end{eqnarray}
\end{corollary}
\pf
By virtue of Theorem~\ref{figo2}, it is enough to prove that, for every~$\Omega^\prime\Subset\Omega$ and any mapping~$G=(g_1,\ldots,g_m)\in\overline{{\mathcal B}(F)}$,
\begin{eqnarray}\label{magnifica2}
\sup_{z\in T_{\Omega^\prime}}
\max_{1\leq\ell\leq m}
\vert
Q_j(g_\ell)(z)-g_\ell(z)
\vert
=
\sup_{z\in T_{\Omega^\prime}}
\max_{1\leq\ell\leq m}
\vert
Q_j(f_\ell)(z)-f_\ell(z)
\vert\,\,,
\end{eqnarray}
or equivalently 
\begin{eqnarray}\label{magnifica3}
\max_{1\leq\ell\leq m}
\sup_{z\in T_{\Omega^\prime}}
\vert
Q_j(g_\ell)(z)-g_\ell(z)
\vert
=
\max_{1\leq\ell\leq m}
\sup_{z\in T_{\Omega^\prime}}
\vert
Q_j(f_\ell)(z)-f_\ell(z)
\vert\,\,.
\end{eqnarray}
So let us fix~$1\leq\ell\leq m$ and let~$(t_r)\subset\mathbb R^n$ be such that~$f_\ell(z+t_r)$ converges uniformly on~$T_{\Omega^\prime}$ to~$g_\ell(z)$, as~$r\to+\infty$. It is easily proved that, for every~$j\in\mathbb N$,
$$
Q_j(g_\ell)(z)
=
\lim_{r\to+\infty}
Q_j(f_{\ell})(z+t_r)
\,,
$$
hence, for every~$z\in T_{\Omega^\prime}$,
\begin{eqnarray*}
\vert
Q_j(g_\ell)(z)-g_\ell(z)
\vert
&=&
\lim_{r\to+\infty}
\vert
Q_j(f_{\ell})(z+t_r)-f_{\ell}(z+t_r)
\vert
\\
&\leq&
\sup_{t\in\mathbb R^n}
\vert
Q_j(f_{\ell})(z+t)-f_{\ell}(z+t)
\vert
\\
&\leq&
\sup_{z\in T_{\Omega^\prime}}
\vert
Q_j(f_{\ell})(z)-f_{\ell}(z)
\vert\,.
\end{eqnarray*}
It follows that
$$
\max_{1\leq\ell\leq m}
\sup_{z\in T_{\Omega^\prime}}
\vert
Q_j(g_\ell)(z)-g_\ell(z)
\vert
\leq
\max_{1\leq\ell\leq m}
\sup_{z\in T_{\Omega^\prime}}
\vert
Q_j(f_{\ell})(z)-f_{\ell}(z)
\vert\,.
$$
Changing the role of~$F$ and~$G$ yields the opposite inequality and the proof is thus complete.\hfill$\square$

\par
\medskip
The following Corollary generalizes the Proposition~\ref{facilis} to the action of~$\Ch\mathbb R^n$ on the entire space~$HAP(T_\Omega,\mathbb C^m)$.
\begin{corollary}\label{2facilis}
let~$F=(f_1,\ldots,f_m)\in HAP(T_\Omega,\mathbb C^m)$ and let~$(\chi_r)\subset\Ch\mathbb R^n$ be a sequence that converges to some~$\chi\in\Ch \mathbb R^n$. Then~$(F_{\chi_r})$ converges to~$F_\chi$ in the topology~$\tau(T_\Omega)$.
\end{corollary}
\pf
Suppose by contradiction the existence of an~$\Omega^\prime\Subset\Omega$ and of an~$\varepsilon>0$ such that, for any~$\nu\in\mathbb N$, there is an integer~$r_\nu>\nu$ and a point~$z_\nu\in T_{\Omega^\prime}$ for which
\begin{eqnarray}\label{continua}
\max_{1\leq\ell\leq m}
\vert
(f_\ell)_\chi(z_\nu)-(f_\ell)_{\chi_{r_\nu}}(z_\nu)
\vert
>\varepsilon\,.
\end{eqnarray}
Then, for every~$\nu$,
\begin{eqnarray*}
\max_{1\leq\ell\leq m}
\vert
(f_\ell)_\chi(z_\nu)-(f_\ell)_{\chi_{r_\nu}}(z_\nu)
\vert
&\leq&
\sup_{z\in T_{\Omega^\prime}}
\max_{1\leq\ell\leq m}
\vert
(f_\ell)_\chi(z)-Q_j(f_\ell)_\chi(z)
\vert
\\
&+&
\sup_{z\in T_{\Omega^\prime}}
\max_{1\leq\ell\leq m}
\vert
Q_j(f_\ell)_\chi(z)-Q_j(f_\ell)_{\chi_{r_\nu}}(z)
\vert
\\
&+&
\sup_{z\in T_{\Omega^\prime}}
\max_{1\leq\ell\leq m}
\vert
Q_j(f_\ell)_{\chi_{r_\nu}}(z)-(f_\ell)_{\chi_{r_\nu}}(z)
\vert\,.
\end{eqnarray*}
Now, by Corollary~\ref{utile1}, there is a~$\tilde\jmath\in\mathbb N$ (depending on~$\Omega^\prime$ and~$\varepsilon$ but not on~$\chi$ or~$\nu$) such that both the first and the third terms in the right-hand side of the above inequality become smaller than~$(\varepsilon/3)$ as soon as~$j>\tilde\jmath\,$. Moreover, if such a~$j$ is fixed, then, by Proposition~\ref{facilis}, there exists an~$\tilde \nu\in\mathbb N$, (depending on~$\Omega^\prime$, on~$\varepsilon$ and on~$j$), such that, when~$\nu>\tilde \nu$, the second term in the right-hand side of the above inequality becomes smaller than~$(\varepsilon/3)$ too. It follows that, when~$j>\tilde\jmath$ and~$\nu>\tilde\nu$, the left-hand side of~(\ref{continua}) is strictly smaller than~$\varepsilon$, whereas it was supposed to be bigger than~$\varepsilon$.\hfill$\square$

\par
\medskip
The following corollary will be used in two important points of the paper.
\begin{corollary}\label{utile2}
Let~$F=(f_1,\ldots,f_m)\in HAP(T_\Omega,\mathbb C)$ and let~$Q_j(F)$ the Bochner-Fejer approximation of~$F$ corresponding to some choice of a basis for~$\lin_{\mathbb Q}\Xi_F$. If~$(\chi_j)\subset\Ch\mathbb R^n$ and~$(\eta_j)\subset \mathbb R^n$ are sequences which converge respectively to~$\chi\in\Ch\mathbb R^n$ and~$\eta\in\mathbb R^n$, then
$$
\lim_{j\to+\infty}
Q_j(F)_{\chi_j}(i\eta_j)
=
F_\chi(i\eta)\,.
$$
\end{corollary}
\pf
Let~$B\Subset\Omega$ be an open ball around~$i\eta$ and let~$j_0$ be such that~$i\eta_j\in B$ when~$j>j_0$. For every such~$j$, we have
\begin{eqnarray*}
\max_{1\leq\ell\leq m}
\vert
Q_j(f_\ell)_{\chi_j}(i\eta_j)-(f_\ell)_\chi(i\eta)
\vert
&\leq&
\sup_{z\in T_B}
\max_{1\leq\ell\leq m}
\vert
Q_j(f_\ell)_{\chi_j}(z)-(f_\ell)_{\chi_j}(z)
\vert
\\
&+&
\sup_{z\in T_B}
\max_{1\leq\ell\leq m}
\vert
(f_\ell)_{\chi_j}(z)-(f_\ell)_\chi(z)
\vert
\\
&+&
\max_{1\leq\ell\leq m}
\vert
(f_\ell)_{\chi_j}(i\eta_j)-(f_\ell)_\chi(i\eta)
\vert
\,.
\end{eqnarray*}
Fix~$\varepsilon>0$. By Corollary~\ref{utile1}, for any~$j$, 
\begin{eqnarray*}
\sup_{z\in T_B}
\max_{1\leq\ell\leq m}
\vert
Q_j(f_\ell)_{\chi_j}(z)-(f_\ell)_{\chi_j}(z)
\vert
=
\sup_{z\in T_B}
\max_{1\leq\ell\leq m}
\vert
Q_j(f_\ell)(z)-f_\ell(z)
\vert\,,
\end{eqnarray*}
so, let~$j_1\in\mathbb N$ be such that
\begin{eqnarray}\label{pezzo1}
\sup_{z\in T_B}
\max_{1\leq\ell\leq m}
\vert
Q_j(f_\ell)(z)-f_\ell(z)
\vert
<(\varepsilon/3)
\,,
\end{eqnarray}
when~$j>j_1$. Moreover, by Corollary~\ref{2facilis}, there is a~$j_2\in\mathbb N$ such that
\begin{eqnarray}\label{pezzo2}
\sup_{z\in T_B}
\max_{1\leq\ell\leq m}
\vert
(f_\ell)_{\chi_j}(z)-(f_\ell)_\chi(z)
\vert
<(\varepsilon/3)
\end{eqnarray}
when~$j>j_2$. Finally, by the continuity of~$F_\chi$ at the point~$i\eta$, there is a~$j_3\in\mathbb N$ such that
\begin{eqnarray}\label{pezzo3}
\max_{1\leq\ell\leq m}
\vert
(f_\ell)_{\chi_j}(i\eta_j)-(f_\ell)_\chi(i\eta)
\vert
<(\varepsilon/3)\,,
\end{eqnarray}
for any~$j>j_3$. Now, if~$j>\max(j_0,j_1,j_2,j_3)$, the relations~(\ref{pezzo1}),~(\ref{pezzo2}) and~(\ref{pezzo3}) imply that
$$
\max_{1\leq\ell\leq m}
\vert
Q_j(f_\ell)_{\chi_j}(i\eta_j)-(f_\ell)_\chi(i\eta)
\vert
<\varepsilon\,.
$$
The Corollary follows by the arbitrary choice of~$\varepsilon$.\hfill$\square$

\section{Amoebae}\label{Amoebae}
We start this section with an easy though important remark.
\begin{lemma}\label{Bochameba}
Let~$F\in HAP(T_\Omega,\mathbb C^m)$ be a regular mapping. Then
$$
{\mathcal F}_F={\mathcal F}_G\,,
$$
for any~$G\in\overline{{\mathcal B}(F)}$.
\end{lemma}
\pf
By~Theorem~\ref{figo2},~${\mathcal O}(F)=\overline{{\mathcal B}(F)}$, so~$\overline{{\mathcal B}(F)}=\overline{{\mathcal B}(G)}$, for any~$G\in \overline{{\mathcal B}(F)}$. It is thus enough to prove the inclusion~${\mathcal F}_G\subseteq{\mathcal F}_F$.

If~$G$ has the form~$G(z)=F(z+t)$, for some~$t\in\mathbb R^n$, the theorem is trivial, so let us suppose
$$
G(z)
=
\lim_{s\to +\infty}
F(z+t_s)
$$
in the topology~$\tau(T_\Omega)$, for some sequence~$(t_s)\subset\mathbb R^n$. By Rouché theorem, the regularity of~$F$ implies that either the zero sets~$V(F)$ and~$V(G)$ are both empty, (in which case the lemma is again trivial), or they are both non empty of codimension~$m$. 
Suppose, by contradiction, there is a point~$y\in {\mathcal F}_G\setminus{\mathcal F}_F$ and let~$\Omega^\prime\Subset\Omega$ be a small open~$n$-dimensional ball around~$y$ disjoint from~${\mathcal F}_F$, so that the tube~$T_{\Omega^\prime}$ is disjoint from~$V(F)$.
Also, let~$S$ an~$m$-dimensional complex affine subspace of~$\mathbb C^n$ passing through~$iy$ such that~$V(G)\cap S$ is discrete. By construction, when~$s$ is big enough, the zero set of~$F(z+t_s)$ cannot intersect~$S$ inside the tube~$T_{\Omega^\prime}$ so, by Rouché theorem, the same has to happen to the limiting mapping~$G$. This is undoubtedly false.\hfill$\square$
\par
\begin{remark}
{\rm
As it will become evident from Remark~3 of this section, the Lemma~\ref{Bochameba} is false without the assumption of regularity for the mapping~$F$.}
\end{remark}

\begin{corollary}
Let~$f\in HAP(T_\Omega,\mathbb C)$, then~${\mathcal F}_f={\mathcal F}_g$, for any~$g\in\overline{{\mathcal B}(f)}$.
\end{corollary}
\pf
A single holomorphic almost periodic function is always regular.\hfill$\square$
\par
\medskip
Lemma~\ref{Bochameba} shows that, for a regular mapping~$F\in HAP(T_\Omega,\mathbb C^m)$, the amoeba~${\mathcal F}_F$ is something that depends only on the closure of the Bochner transform of~$F$ and not on the mapping~$F$ itself. Now, for a regular~$F$,
$$
{\mathcal F}_F
=
\bigcup_{G\in\overline{{\mathcal B}(F)}}
\overline{\im V(G)}\,,
$$
and since~$\im V(F)=\im V(G)$, for any~$G\in {\mathcal B}(G)$, one may ask if the points from~${\mathcal F}_F\setminus\im V(F)$ are zeros of some mapping~$G\in\overline{{\mathcal B}(F)}\setminus{\mathcal B}(F)$. Equivalently, one may ask if, for a regular~$F$, the following equalities hold true
$$
{\mathcal F}_F
=
\bigcup_{G\in\overline{{\mathcal B}(F)}}
\im V(G)
=
\bigcup_{\chi\in\Ch\mathbb R^n}
\im V(F_\chi)\,.
$$
Observe that, if the answer to this question is yes, the above equalities would provide a representation of~${\mathcal F}_F$ that is somehow more effective than its definition. In fact, though these equalities require the union of an infinite family, they do not involve any topological closure. In order to give a positive answer to this question, we introduce the following auxiliary notion of amoeba.
\par\medskip
Let~$F\in HAP(T_\Omega,\mathbb C^m)$ and~$\mathbb G$ be an additive subgroup of~$\mathbb R^n$ containing the spectra of the components of~$F$. The ${\mathcal Y}${\sl -amoeba}~${\mathcal Y}_F$ of~$F$ is the union over~$\Ch\mathbb G$ of the imaginary parts of the zeros of the perturbed mappings, namely
$$
{\mathcal Y}_F
=
\bigcup_{\chi\in\Ch\mathbb G}
\im V(F_\chi)\,.
$$
\par
\begin{remark}
{\rm
By the injectiveness of the group~$\mathbb S^1$, (see~Remark~2 of section~3), the~${\mathcal Y}$-amoeba of~$F$ does not depend on the choice of the group~$\mathbb G$, so~${\mathcal Y}_F$ is a well defined object that can be described, depending on the particular situation, by using the characters of the most convenient group~$\mathbb G$.}
\end{remark}
\par

\begin{theorem}\label{pratico}
Let~$F\in HAP(T_\Omega,\mathbb C^m)$. Then the ${\mathcal Y}$-amoeba of~$F$ is a closed subset such that
\begin{eqnarray}\label{prima}
i{\mathcal Y}_F
=
i\Omega
\cap
\bigcup_{\chi\in\Ch\mathbb R^n}
V(F_\chi)\,.
\end{eqnarray}
\end{theorem}
\pf
We first show the relation~(\ref{prima}). 
Let~$\eta\in{\mathcal Y}_F$, so there exists~$\xi\in\mathbb R^n$ and a~$\chi\in\Ch\mathbb R^n$ such that
$$
0=
F_\chi(\xi+i\eta)
=
\lim_{j\to+\infty}
F_{j,\chi}(\xi+i\eta)\,,
$$
where~$(F_j)\subset\Ex(T_\Omega,\mathbb C^m)$ is the Bochner-Fejer approximation of~$F$ corresponding to some choice of a basis for~$\lin_{\mathbb Q}\Xi_F$. If~$\chi^\prime$ denotes the character of~$\mathbb R^n$ given by~$e^{i\langle \xi,\cdot\rangle}$, we get
$$
0=
F_{\chi\chi^\prime}(i\eta)
=
\lim_{j\to+\infty}
(F_{j,\chi})_{\chi^\prime}(i\eta)
\,,
$$
whence~$i\eta\in i\Omega\cap V(F_{\chi\chi^\prime})$. The other inclusion is trivial.
\par
Let~$(\eta_j)\subset {\mathcal Y}_F$ a sequence converging to~$\eta\in \mathbb R^n$. For every~$j\in\mathbb N$, there exists a~$\chi_j\in\Ch\Xi_F$ such that~$i\eta_j\in V(F_{\chi_j})$. By compactness of~$\Ch\mathbb R^n$, the sequence~$(\chi_j)$ has a convergent subsequence~$(\chi_{j_k})$. Let~$\chi$ be the limiting character, then, by Corollary~\ref{utile2}
$$
F_\chi(i\eta)
=
\lim_{k\to+\infty}
F_{\chi_{j_k}}
(i\eta_{j_k})
=0\,,
$$
whence~$\eta\in{\mathcal Y}_F$.\hfill$\square$

\par
\begin{remark}\label{inter}
{\rm 
For every~$F\in HAP(T_\Omega,\mathbb C^m)$ and any~$\chi\in\Ch\mathbb R^n$, the inclusions 
$$
\im V(F_\chi)
\subseteq
{\mathcal F}_{F_\chi}
\subseteq
{\mathcal Y}_F
$$
hold true, nevertheless they generally hold strictly (cf.~\cite{sil06}). Moreover~Theorem~\ref{pratico} implies the equality
$$
{\mathcal Y}_F
=
\bigcup_{\chi\in\Ch\mathbb R^n}
{\mathcal F}_{F_\chi}\,,
$$
for any~$F\in HAP(T_\Omega,\mathbb C^m)$.}
\end{remark}
\begin{corollary}\label{iorfigo}
Let~$1\leq m\leq n$ and~$F\in HAP(T_\Omega,\mathbb C^m)$ be a regular mapping. Then~${\mathcal Y}_F={\mathcal F}_F$.
\end{corollary}
\pf
By~Theorem~\ref{figo2} and~Lemma~\ref{Bochameba},
$$
{\mathcal Y}_F\subseteq\bigcup_{\chi\in\Ch\mathbb R^n}\overline{\im V(F)}={\mathcal F}_F\,.
$$
On the other hand, by~Theorem~\ref{pratico},~${\mathcal Y}_F$ is closed and as~$\im V(F)\subseteq{\mathcal Y}_F$, we get~${\mathcal F}_F\subseteq{\mathcal Y}_F$.\hfill$\square$
\par
\medskip
The last two lemmas of this section deal with amoebae of holomorphic almost periodic mappings whose components are entire functions.
\begin{lemma}\label{cambiovar}
Let~$F\in HAP(\mathbb C^n,\mathbb C^m)$ and~$\varphi$ a~$\mathbb C$-linear automorphism of~$\mathbb C^n$ such that~$\varphi(\mathbb R^n)=\mathbb R^n$. Then
\begin{eqnarray}\label{i}
\im[V(F\circ\varphi)]=\im[\varphi^{-1}(V(F))]=\varphi^{-1}(\im V(F))
\end{eqnarray}
and
\begin{eqnarray}\label{ii}
{\mathcal Y}_F
=
\varphi^a({\mathcal Y}_{F\circ\varphi^a})\,,
\end{eqnarray}
where~$\varphi^a$ is the adjoint of~$\varphi$ with respect to the standard hermitian form on~$\mathbb C^n$.
\end{lemma}
\pf
The proof of~(\ref{i}) is straightforward, so let us show the equality~(\ref{ii}).
If~$F\in\Ex(\mathbb C^n,\mathbb C^m)$, the proof is formally the same as the one given in the Lemma~3.2 of the paper~\cite{sil06}. Let~$1\leq\ell\leq m$ and let~$f_\ell$ be the~$\ell$-th component of~$F$. We have
$$
f_\ell(z)
=
\sum_{\lambda\in\Sp f_\ell}
a(\lambda,f_\ell)\,
e^{i\langle z,\lambda\rangle}\,,
$$
so, for any~$\chi\in\Ch\mathbb R^n$,
\begin{eqnarray*}
(f_\ell)_\chi\circ\varphi^a(z)
&=&
\sum_{\lambda\in\Sp f_\ell}
a(\lambda,(f_\ell)_\chi)\,
e^{i\langle\varphi^a(z),\lambda\rangle}
=
\sum_{\lambda\in\Sp f_\ell}
a(\lambda,(f_\ell)_\chi)\,
e^{i\overline{\langle\lambda,\varphi^a(z)\rangle}}
\\
&=&
\sum_{\lambda\in\Sp f_\ell}
a(\lambda,(f_\ell)_\chi)\,
e^{i\overline{\langle\varphi(\lambda),z\rangle}}
=
\sum_{\lambda\in\Sp f_\ell}
a(\lambda,(f_\ell)_\chi)\,
e^{i\langle z,\varphi(\lambda)\rangle}
\\
&=&
\sum_{\lambda\in\Sp f_\ell}
a(\lambda,f_\ell)\,
\chi(\lambda)\,
e^{i\langle z,\varphi(\lambda)\rangle}
\,,
\end{eqnarray*}
which implies that~$\Sp(f_\ell\circ\varphi^a)=\varphi(\Sp f_\ell)$ and
$$
a(\varphi(\lambda),f_\ell\circ\varphi^a)=a(\lambda,f_\ell)\,,
$$
for every~$\lambda\in\Sp f_\ell$. 
On the other hand,~$\chi\circ\varphi^{-1}\in\Ch\mathbb R^n$, for any~$\chi\in\Ch\mathbb R^n$, and
\begin{eqnarray*}
(f_\ell\circ\varphi^a)_{\chi\circ\varphi^{-1}}(z)
&=&
\sum_{\varphi(\lambda)\in\varphi(\Sp f_\ell)}
a(\varphi(\lambda),f_\ell\circ\varphi^a)
\,
\chi\circ\varphi^{-1}(\varphi(\lambda))
\,
e^{i\langle z,\varphi(\lambda)\rangle}
\\
&=&
\sum_{\lambda\in\Sp f_\ell}
a(\lambda,f_\ell)
\,
\chi(\lambda)
\,
e^{i\langle z,\varphi(\lambda)\rangle}\,,
\end{eqnarray*}
for every~$1\leq\ell\leq m$, so
$$
F_\chi\circ\varphi^a
=
(F\circ\varphi^a)_{\chi\circ\varphi^{-1}}\,,
$$
and, thanks to~(\ref{i}),
$$
\im V(F_\chi)
=
\varphi^a\big(\im V(F\circ\varphi^a)_{\chi\circ\varphi^{-1}}\big)\,.
$$
Taking the union over~$\chi\in\Ch\mathbb R^n$ on both sides of the above equality, we get~${\mathcal Y}_F=\varphi^a({\mathcal Y}_{F\circ\varphi^a})$. This proves the lemma in the case~$F\in\Ex(\mathbb C^n,\mathbb C^m)$. In the general case, let us consider the Bochner-Fejer approximation~$Q_j(F)$ of~$F$ corresponding to some choice of a basis of~$\lin_{\mathbb Q}\Xi_F$. For any~$\chi\in\Ch\Xi_F$, thanks to the previous discussion, we have
$$
F_\chi\circ\varphi^a
=
\lim_{j\to+\infty}
Q_j(F)_\chi\circ\varphi^a
=
\lim_{j\to+\infty}
(Q_j(F)\circ\varphi^a)_{\chi\circ\varphi^{-1}}
=
(F\circ\varphi^a)_{\chi\circ\varphi^{-1}}\,,
$$
so~$\im V(F_\chi\circ\varphi^a)=\im V((F\circ\varphi^a)_{\chi\circ\varphi^{-1}})$. Using~(\ref{i}) and taking the union over~$\chi\in\Ch\mathbb R^n$ completes the proof.\hfill$\square$
\par
\begin{lemma}\label{rat}
If~$F=(f_1,\ldots,f_m)\in HAP(\mathbb C^n,\mathbb C^m)$, let~$\lin_\mathbb R\Xi_F$ be the real linear subspace of~$\,\mathbb R^n$ generated by~$\Xi_F$. If~$\,\Xi_F$ is free of finite rank, then
$$
{\mathcal Y}_F=\im V(F)={\mathcal F}_F\,,
$$
as soon as the rank of~$\,\Xi_F$ equals the dimension of~$\,\lin_\mathbb R \Xi_F$. In particular the above equality holds true when~$\Xi_F$ is a finitely generated subgroup of~$\mathbb Q^n$.
\end{lemma}
\pf
We first suppose that~$\Xi_F$ has the following form
$$
\Xi_F
=
\{(p_1,\ldots,p_s,0,\ldots,0)\in\mathbb R^n\mid p_1,\ldots,p_s\in\mathbb Z\}\,,
$$
where~$s=\rk\Xi_F\leq n$. In this case we may represent the amoeba of~$F$ by means of the chracters of~$\mathbb Z^n$. Each such character is uniquely determined by its values on the canonical base of~$\mathbb Z^n$. If~$\chi\in\Ch\mathbb Z^n$ is the character sending~$p\in\mathbb Z^n$ to~$e^{i(p_1\vartheta_1+\cdots+p_n\vartheta_n)}$, then, for every~$\lambda=(\lambda_1,\ldots,\lambda_s,0,\ldots,0)$ in~$\Xi_F$ we have
\begin{eqnarray*}
\chi(\lambda)\,e^{i\langle z,\lambda\rangle}
&=&
e^{i(\lambda_1\vartheta_1+\cdots+\lambda_s\vartheta_s)}
e^{i(z_1\lambda_1+\cdots+z_s\lambda_s)}
\\
&=&
e^{i(z_1+\vartheta_1)\lambda_1+\cdots+(z_s+\vartheta_s)\lambda_s}
\\
&=&
e^{i\langle z+\vartheta,\lambda\rangle}\,.
\end{eqnarray*}
If~$Q_j(F)$ is the Bochner-Fejer approximation of~$F$ that corresponds to the canonical basis of~${\mathbb Q}^s=\lin_{\mathbb Q}\Xi_F$, it follows that
$$
F_\chi(z)
=
\lim_{j\to+\infty}
Q_j(F)_\chi(z)
=
\lim_{j\to+\infty}
Q_j(F)(z+\vartheta)
=
F(z+\vartheta)\,,
$$
so~$\im V(F_\chi)=\im V(F)$. Since~$\chi\in\Ch\Xi_F$ was arbitrary, taking the union over~$\chi\in\Ch\Xi_F$ yelds~${\mathcal Y}_F=\im V(F)$. 
\par
Let us now prove the general case. Let~$s=\rk\Xi_F=\dim_\mathbb R\lin_\mathbb R\Xi_F$ and let~$\{\omega_1,\ldots,\omega_s\}$ a system of free generators of~$\Xi_F$. This elements of~$\mathbb R^n$ are~$\mathbb R$-linearly independent otherwise the~$\mathbb R$-linear subspace they span, namely~$\lin_\mathbb R\Xi_F$, would have a dimension strictly smaller than~$s$. Let us complete the system~$\{\omega_1,\ldots,\omega_s\}$ to a basis~$\{\omega_1,\ldots,\omega_n\}$ of~$\mathbb R^n$ and let~$\varphi$ be the~$\mathbb C$-linear automorphism of~$\mathbb C^n$ preserving~$\mathbb R^n$ and represented, with respect to the canonical basis, by the inverse of the matrix having the vectors~$\omega_1,\ldots,\omega_n$ as columns. Up to an odd permutation of these columns, we may also assume~$\det\varphi>0$. By construction, we have
$$
\varphi(\Xi_F)
=
\{(p_1,\ldots,p_s,0,\ldots,0)\in\mathbb R^n\mid p_1,\ldots,p_s\in\mathbb Z\}\,,
$$
so, by the first part of the proof, we know that~${\mathcal Y}_{F\circ\varphi}=\im V(F\circ\varphi)$. Lemma~\ref{cambiovar} implies that
$$
{\mathcal Y}_F
=
\varphi({\mathcal Y}_{F\circ\varphi})
=
\varphi(\im V(F\circ\varphi))
=
\varphi(\varphi^{-1}(\im V(F)))
=
\im V(F)\,.
$$
Finally, if all the spectra of the components of~$F$ are contained in~$\mathbb Q^n$, we may complete the system~$\{\omega_1,\ldots,\omega_s\}$ of free generators of~$\Xi_F$ to a basis~$\{\omega_1,\ldots,\omega_n\}$ of~$\mathbb Q^n$. Now construct, as before, a~$\mathbb C$-linear authomorphism~$\varphi$ of~$\mathbb C^n$. It follows that~$\varphi$ preserves~$\mathbb Q^n$ and
$$
\varphi(\lin_\mathbb Q\Xi_F)
=
\{(p_1,\ldots,p_s,0,\ldots,0)\in\mathbb R^n\mid p_1,\ldots,p_s\in\mathbb Q\}\,,
$$
hence
$$
\dim_\mathbb R\lin_\mathbb R\Xi_F
=
\dim_\mathbb R\varphi(\lin_\mathbb R\Xi_F)
=
\dim_\mathbb Q\varphi(\lin_\mathbb Q\Xi_F)
=
s\,,
$$
so the conclusion follows from the second part of the proof.\hfill$\square$

\par
\begin{remark}
{\rm 
The condition stated in~Lemma~\ref{rat} is sufficient but not necessary. In fact, let~$F_1,F_2\in HAP(\mathbb C^n,\mathbb C)$ be such that~$\Xi_{F_1},\Xi_{F_2}$ are free with finite ranks,~$\lin_\mathbb R\Xi_{F_1}=\lin_\mathbb R\Xi_{F_2}$ and~$\rk(\Xi_{F_1}+\Xi_{F_2})=\rk\Xi_{F_1}+\rk\Xi_{F_2}$. Then the holomorphic almost periodic function~$F= F_1 F_2$ does not satisfy the hypothesis of the~Lemma~\ref{rat}, nevertheless~$V(F_\chi)=V((F_1)_\chi)\cup V((F_2)_\chi)$, for any character~$\chi\in\Ch\mathbb R^n$, hence~${\mathcal Y}_F=\im V(F)={\mathcal F}_F$.
\par
Remark also that~Lemma~\ref{rat} is not true if~$\Xi_F$ is not free or if it is free with infinite rank. In fact, there exist examples of functions~$F\in HAP(\mathbb C,\mathbb C)$, with~$\Xi_F=\mathbb Q$ or with~$\Xi_F$ equal to the additive subgroup of~$\mathbb R$ consisting of all the real algebraic numbers, for which the amoeba~${\mathcal F}_F$ equals~$\mathbb R$, while~$\im V(F)$ is only countable.
}
\end{remark}

\section{Amoeba complements}\label{Amoebacompl}

We finally turn our attention toward Henriques convexity properties of amoeba complements. 
\par
\begin{lemma}\label{buonino}
Let~$F=(f_1,\ldots,f_{m+1})\in\Ex(\mathbb C^n,\mathbb C^{m+1})$, with~$m+1\leq n$, such that~$\Xi_F\subseteq\mathbb Q^n$. If~$s=\codim V(F)$, then the amoeba complement~$\mathbb R^n\setminus{\mathcal F}_F$ is an~$(s-1)$-convex subset; in particular it is~$m$-convex.
\end{lemma}
\pf
We may suppose~$\varnothing\neq{\mathcal F}_F\neq\mathbb R^n$ and~$0\leq m<n-1$, otherwise the theorem is trivial. Moreover, if~$m=0$ the theorem specializes to a well-known result. Indeed, if~$m=0$, each connected component~$X$ of~$\mathbb R^n\setminus{\mathcal F}_F$ is the base of a proper, tubular, connected component~$T_X$ of the subset~$\mathbb R^n+i(\mathbb R^n\setminus{\mathcal F}_F)\subset\mathbb C^n$, so~$T_X$ is an open and pseudoconvex tube in~$\mathbb C^n$, which implies the convexity of~$X$. We then suppose~$m\geq 1$.
\par
If~$F\in\Ex(\mathbb C^n,\mathbb C^m)$ and~$\Xi_F\subseteq\mathbb Z^n$, the theorem reduces to~Theorem~\ref{henr} proved by Henriques. In fact, let
$$
P=(p_1,,\ldots,p_{m+1}):(\mathbb C^*)^n\rightarrow\mathbb C^n
$$
be the mapping given, for every~$1\leq\ell\leq (m+1)$ and every~$\zeta\in(\mathbb C^*)^n$, by
$$
p_\ell(\zeta)
=
\sum_{\lambda\in\Sp f_\ell}
a(\lambda,f_\ell)\,
\zeta^{-\lambda}\,,
$$
and let~$\psi:\mathbb C^n\longrightarrow(\mathbb C^*)^n$ be the mapping given, for every $z\in\mathbb C^n$, by
$$
\psi(z_1,\ldots,z_n)=(e^{-iz_1},\ldots,e^{-iz_n})\,.
$$
It is easily seen that~$\psi(V(F))$ equals the zero set~$V$ of~$P$ in the torus~$(\mathbb C^*)^n$, moreover, the mapping~$\psi$ has discrete fibers, so~$\codim V=s$ and by Theorem~\ref{henr} the complementary set to the polynomial amoeba~${\mathcal A}_V$ is~$(s-1)$-convex. Since~${\mathcal F}_F={\mathcal A}_V$, it follows that~$\mathbb R^n\setminus{\mathcal F}_F$ is~$(s-1)$-convex too. As~$m\geq (s-1)$, the set~$\mathbb R^n\setminus{\mathcal F}_F$ is also~$m$-convex.
\par
Finally, if~$F\in\Ex(\mathbb C^n,\mathbb C^m)$ and~$\Sp f_1,\ldots,\Sp f_{m+1}\subset\mathbb Q^n$, let~$\varphi$ be a~$\mathbb C$-linear authomorphism of~$\mathbb C^n$ such that~$\varphi(\mathbb R^n)=\mathbb R^n$,~$\det\varphi>0$ and~$\varphi(\Xi_F)\subset\mathbb Z^n$. By~Lemma~\ref{cambiovar} and Lemma~\ref{rat} we get
$$
\mathbb R^n\setminus{\mathcal F}_F
=
\varphi^a(\mathbb R^n\setminus{\mathcal F}_{F\circ\varphi^a})\,.
$$
Now, the spectra of the components of~$F\circ\varphi^a$ are contained in~$\mathbb Z^n$ and~$m$-convexity is preserved by~$\mathbb R$-linear authomorphisms of~$\mathbb R^n$ with positive determinant, so the previous case makes the theorem hold true in the present one too.\hfill$\square$

\par
\begin{theorem}\label{buono}
Let~$F=(f_1,\ldots,f_{m+1})\in \Ex(\mathbb C^n,\mathbb C^{m+1})$ with~$m+1\leq n$ and let~$\Omega\subseteq\mathbb R^n$ be a non-empty open and convex subset. If~$F$ is regular on~$T_\Omega$, the set~$\Omega\setminus{\mathcal F}_F\subseteq\mathbb R^n$ is an~$m$-convex subset.
\end{theorem}
\pf
As in~Lemma~\ref{buonino} ,we suppose~$\varnothing\neq{\mathcal F}_F\neq\Omega$ and~$1\leq m<n-1$.
We use a density argument involving holomorphic almost periodic mappings whose~$\mathbb Z$-modules are included in~$\mathbb Q^n$.
We split the proof into several steps. 
\par\smallskip
\noindent
\textit{Step~1.}
Let~$\omega_1,\ldots,\omega_r$ a system of free generators of the~$\mathbb Z$-module~$\Xi_F$. As in the proof of~Lemma~\ref{figo}, for any~$1\leq k\leq (m+1)$, we have
\begin{eqnarray}
f_k(z)
=
\sum_{\mu\in\mathbb Z^r}
a_{k,\mu}
\prod_{\ell=1}^r
e^{i\mu_\ell\langle z,\omega_\ell\rangle}\,,
\end{eqnarray}
where all but a finite number of~$a_{k,\mu}$'s are different from zero. For every character~$\chi\in\Ch\Xi_F$, we correspondingly get the expression
\begin{eqnarray}
(f_k)_\chi(z)
=
\sum_{\mu\in\mathbb Z^r}
a_{k,\mu}\,
\chi\biggl(\sum_{\ell=1}^r\mu_\ell\,\omega_\ell\biggr)
\prod_{\ell=1}^r
e^{i\mu_\ell\langle z,\omega_\ell\rangle}\,.
\end{eqnarray}
Now, for every~$1\leq\ell\leq r$, let~$(\omega_{\ell,j})_{j\in\mathbb N}\subset\mathbb Q^n$ a sequence converging to~$\omega_\ell$ and, for any~$j\in\mathbb N$ and any~$\chi\in\Ch\Xi_F$, let~$(F_\chi)_j\in\Ex(\mathbb C^n,\mathbb C^{m+1})$ the mapping whose~$k$-th component~$((f_k)_\chi)_j$ is equal to
\begin{eqnarray}
((f_k)_\chi)_j(z)
=
\sum_{\mu\in\mathbb Z^r}
a_{k,\mu}\,
\chi\biggl(\sum_{\ell=1}^r\mu_\ell\,\omega_\ell\biggr)
\prod_{\ell=1}^r
e^{i\mu_\ell\langle z,\omega_{\ell,j}\rangle}\,.
\end{eqnarray}
We thus obtain, for any~$\chi\in\Ch\Xi_F$, a sequence~$(F_\chi)_j$ of mappings, with spectra contained in~$\mathbb Q^n$, converging to~$F_\chi$ uniformly on compacta, as~$j\to+\infty$.
\par
Moreover, by~Lemma~\ref{buonino}, for any~$j$ and every~$\chi$, the set~$\mathbb R^n\setminus\im V((F_\chi)_j)$ is~$m$-convex. 
\par\smallskip
\noindent
\textit{Step~2.}
Suppose, by contradiction, that~$\Omega\setminus {\mathcal F}_F$ is not~$m$-convex. Since we have~$\varnothing\neq{\mathcal F}_F\neq\Omega$ and~$m<n-1$, there exists at least an oriented affine subspace~$S\subset\mathbb R^n$ of dimension~$(m+1)$ such that~$\varnothing\neq S\setminus{\mathcal F}_F\neq S\cap\Omega$. Among such subspaces, we may find one, say~$S$, for which there exists at least a polyhedral chain~$c\in {\mathcal C}^\Delta_m (S\cap\Omega\setminus {\mathcal F}_F)$ that represents a non trivial non negative homology class~$[c]\in H_m^+(S\cap\Omega\setminus{\mathcal F}_F)$ becoming homologically trivial in~$\Omega\setminus{\mathcal F}_F$. We will show that the existence of such a chain leads to a contradiction. 
\par
Thanks to Lemma~2.7 of~\cite{hen04}, the chain~$c$ bounds a unique~$C\in{\mathcal C}^\Delta_{m+1}(S)$; its support~$\Supp C$ is not a subset of~$S\cap\Omega\setminus {\mathcal F}_F$. By~Corollary~\ref{iorfigo} this means that there is a~$\chi_o\in\Ch\mathbb R^n$ such that~$\Supp C$ is not a subset of~$S\cap\Omega\setminus \im V(F_{\chi_o})$. Moreover,~$c$ is the boundary of some~$D\in{\mathcal C}^\Delta_{m+1}(\Omega\setminus{\mathcal F}_F)$.
\par\smallskip
\noindent
\textit{Step~3.}
We show that there is a~$k_1\in\mathbb N$ such that, for any~$\chi\in\Ch\Xi_F$, the chain~$c$ is homologically trivial in~$\Omega\setminus{\mathcal F}_{(F_\chi)_j}$ as soon as~$j\geq k_1$. In order to do so, let us prove that there is a~$k_1\in\mathbb N$ such that, for every~$\chi\in\Ch\Xi_F$,
\begin{eqnarray}\label{inclusionechi}
A=\Supp c\cup\Supp D\subset\Omega\setminus\im V((F_\chi)_j)\,,
\end{eqnarray}
as soon as~$j\geq k_1$. We start by remarking that, as~$A$ is compact, there is a finite number of closed balls~$B_1,\ldots,B_t$ such that
\begin{eqnarray}\label{facile}
A\subset B_1\cup\ldots\cup B_t\subset\Omega\setminus{\mathcal F}_F\,.
\end{eqnarray}
If we show that, for any~$1\leq s\leq t$, there is a~$\nu_s\in\mathbb N$ such that, for any~$\chi$,
$$
B_s\subset\Omega\setminus\im V((F_\chi)_j)\,,
$$
when~$j\geq\nu_s$, it will be enough to take~$k_1=\max\{\nu_1,\ldots,\nu_t\}$. So let us prove it by contradiction and let us suppose that, for some~$1\leq s\leq t$, there is an increasing sequence of integer~$j_q$ and a sequence of characters~$\chi_q$, such that~$B_s$ intersects~$\im V(F_{\chi_q})_{j_q}$ for every~$q\in\mathbb N$. For any~$q$, let~$\eta_q$ be a point in that intersection and let~$\xi_q\in\mathbb R^n$ be such that
\begin{eqnarray}\label{equ}
(F_{\chi_q})_{j_q}(\xi_q+i\eta_q)=0\,.
\end{eqnarray}
Now, for every~$1\le k\leq (m+1)$,
\begin{eqnarray*}
((f_k)_{\chi_q})_{j_q}(\xi_q+i\eta_q)
&=&
\sum_{\mu\in\mathbb Z^r}
a_{k,\mu}\,
\chi_q\biggl(\sum_{\ell=1}^r\mu_\ell\,\omega_\ell\biggr)
\prod_{\ell=1}^r
e^{i\mu_\ell\langle \xi_q+i\eta_q,\omega_{\ell,j_q}\rangle}
\\
&=&
\sum_{\mu\in\mathbb Z^r}
a_{k,\mu}\,
\tilde\chi_q\biggl(\sum_{\ell=1}^r\mu_\ell\,\omega_\ell\biggr)
\prod_{\ell=1}^r
e^{i\mu_\ell\langle i\eta_q,\omega_{\ell,j_q}\rangle}\,,
\end{eqnarray*}
where, for any~$q\in\mathbb N$,~$\tilde\chi_q=\chi_q\kappa_q$ and~$\kappa_q$ is the only character of~$\Xi_F$ such that, for any~$1\leq \ell\leq r$,~$\kappa_q(\omega_\ell)=\exp{i\langle\xi_q,\omega_{\ell,j_q}\rangle}$. Consequently equation~(\ref{equ}) becomes
\begin{equation}\label{equa}
(F_{\tilde\chi_q})_{j_q}(i\eta_q)=0\,.
\end{equation}
By the compactness of~$B_s$ and~$\Ch\Xi_F$, the sequences~$(\eta_q)$ and~$(\tilde\chi_q)$ admit subsequences~$(\eta_{q_p})$ and~$(\tilde\chi_{q_p})$ converging respectively to some~$\eta\in B_s$ and some~$\tilde\chi\in\Ch\Xi_F$ as~$p\to+\infty$, so
$$
F_{\tilde\chi}(i\eta)
=
\lim_{p\to+\infty}(F_{\tilde\chi_{q_p}})_{j_{q_p}}(i\eta_{q_p})
\,,
$$
and~(\ref{equa}) yields~$F_{\tilde\chi}(i\eta)=0$, which is impossible since~$\eta\in B_s\subset\Omega\setminus{\mathcal F}_F$.
\par
\smallskip\noindent
\textit{Step~4.}
If~$j$ is bigger than~$k_1$,~$y$ is a point of~$S\cap\im V(F_{\chi_o})_j$ and~$\upsilon_y$ is the standard generator of the de Rham cohomology group~$H_{dR}^m(S\setminus\{y\})$, it is easily proved that
$$
\int_c\upsilon_y>0\qquad{\hbox{\rm if~$y\in\Supp C$,}}
\qquad
\hbox{\rm whereas}
\qquad
\int_c\upsilon_y=0\qquad{\hbox{\rm if~$y\not\in\Supp C$}}\,.
$$
Stokes theorem implies that, for any~$j\geq k_1$, the chain~$c$ represents a non negative homology class in the group~$H_m(S\setminus\im V((F_{\chi_o})_j))$, so if~$j$ is bigger than~$k_1$, the~$m$-convexity of~$\mathbb R^n\setminus\im V((F_{\chi_o})_j)$  implies that~$c$ is homologically trivial in~$S\setminus\im V((F_{\chi_o})_j)$, or equivalently that
\begin{eqnarray}\label{assurdo}
\Supp C\subset S\setminus\im V((F_{\chi_o})_j)\,,
\end{eqnarray}
for any~$j\geq k_1$. Moreover, since~$\Omega$ is convex and~$\Supp c\subset\Omega$, it follows that~$\Omega$ contains the convex hull of the set~$\Supp c$, which of course admits~$\Supp C$ as a subset. In particular, this implies that
\begin{eqnarray}\label{Assurdo}
\Supp C\subset S\cap\Omega\setminus\im V((F_{\chi_o})_j)\,,
\end{eqnarray}
as soon as~$j\geq k_1$.
\par
\smallskip
\noindent
\textit{Step~5.}
Let~$y\in S\cap\Omega\cap\im V(F_{\chi_o})$ be a point in the relative interior of~$\Supp C$ and let~$W$ be a neighborhood of~$y$ such that~$W\cap S\subset\Supp C$. If~$E_S$ is the linear subspace parallel to~$S$ and~$x\in\mathbb R^n$ is such that~$z= x+iy\in V(F_{\chi_o})$, let us consider the complex affine~$(m+1)$-dimensional subspace
$$
U=(x+E_S)+iS\,.
$$
Since~$S\setminus{\mathcal F}_F\neq\varnothing$ it is clear that~$U\not\subset V(F_{\chi_o})$, and as the analytic subset~$V(F_{\chi_o})\cap T_\Omega$ is~$(m+1)$-codimensional, it follows that~$V(F_{\chi_o})\cap T_\Omega$ intersects~$U$ along a discrete analytic subset of~$\mathbb C^n$. 
Let~$B_z$ an open ball about~$z$ not containing other points of~$U\cap V(F_{\chi_o})\cap T_\Omega$. By Rouché theorem, there is a~$k_2\in\mathbb N$, (depending on~$B_z$), such that, for any~$j\geq k_2$, the analytic subset~$U\cap V((F_{\chi_o})_j)\cap T_\Omega$ admits a single point, say~$z_j= x_j+iy_j$, inside the ball~$B_z$. Now, since the sequence~$(z_j)$ converges to~$z$, there is a~$k_3\geq k_2$ such that, for every~$j$ bigger than~$k_3$, the point~$y_j$ belongs to the set~$W\cap S$, which by construction is included in~$\Supp C$. 

This fact yields a contradiction with the inclusion~(\ref{Assurdo}) as soon as~$j$ is bigger than~$k_4=\max(k_1,k_3)$, in fact, in this case, the point~$y_j\in S\cap\Omega\cap \im V(F_{\chi_o})_j$ should also belong to~$S\cap\Omega\setminus \im V(F_{\chi_o})_j$.\hfill$\square$

\par

\begin{corollary}\label{buonis}
Let~$\Omega\subseteq\mathbb R^n$ be a non-empty, open, convex subset and let also~$F=(f_1,\ldots,f_{m+1})\in HAP(T_\Omega,\mathbb C^{m+1})$, with~$m+1\leq n$. If~$F$ is regular, then~$\Omega\setminus{\mathcal F}_F$ is an~$m$-convex subset.
\end{corollary}
\pf
As in~Lemma~\ref{buonino} and~Theorem~\ref{buono}, we may suppose~$\varnothing\neq{\mathcal F}_F\neq\mathbb R^n$ and~$1\leq m<n-1$.
Here we use again a density argument, but now we pass through Bochner-Fejer theorem. The structure of the proof is the same as that of~Theorem~\ref{buono}.
\par
\smallskip
\noindent
\textit{Step~1.}
Let~$F_j=Q_j(F)$ be the Bochner-Fejer approximation of~$F$ corresponding to some choice of a basis of~$\lin_{\mathbb Q}\Xi_F$. Observing that, for any~$\chi\in\Ch\mathbb R^n$ and any~$j$, we have~$(Q_j(F))_\chi=Q_j(F_\chi)$, we get that, for any~$\chi$,~$(F_j)_\chi$ is the Bochner-Fejer approximation of~$F_\chi$, so that we may write~$(F_j)_\chi=(F_\chi)_j$, for every~$j$ and~$\chi$.
\par
By virtue of Lemma~\ref{ras}, for every non-empty, open, convex and relatively compact subset~$\Omega^\prime\subseteq\Omega$, there exists a~$k_0\in\mathbb N$, (depending on~$\Omega^\prime$), such that the mapping~$F_j$ is regular on~$T_{\Omega^\prime}$, for every~$j\geq k_0$. If~$j\geq k_0$, Theorem~\ref{buono}, together with Lemma~\ref{Bochameba} and Theorem~\ref{figo2}, applied to the mapping~$F_j\in HAP(T_{\Omega^\prime},\mathbb C^{m+1})$, imply the~$m$-convexity of the set~$\Omega^\prime\setminus\overline{\im V((F_\chi)_j)}$, for every~$\chi\in\Ch\mathbb R^n$.
\par
\smallskip
\noindent
\textit{Step~2.}
We now suppose, by contradiction, that~$\Omega\setminus{\mathcal F}_F$ is not an~$m$-convex set, so we may {\it verbatim\/} repeat the considerations of Step~2 in the proof of~Theorem~\ref{buono}.
\par
\smallskip
\noindent
\textit{Step~3.}
Here we show that there is a~$k_1\in\mathbb N$ such that the chain~$c$ is homologically trivial in~$\Omega\setminus{\mathcal F}_{(F_{\chi_o})_j}$ as soon as~$j\geq k_1$. In order to do so, let us prove that there is a~$k_1\in\mathbb N$ such that 
\begin{eqnarray}\label{inclusionechio}
A=\Supp c\cup\Supp D\subset\Omega\setminus\overline{\im V((F_{\chi_o})_j)}\,,
\end{eqnarray}
as soon as~$j\geq k_1$. We start by remarking that, as~$A$ is compact, there is a finite number of closed balls~$B_1,\ldots,B_t$ such that
\begin{eqnarray}\label{facile2}
A\subset B_1\cup\ldots\cup B_t\subset\Omega\setminus{\mathcal F}_F\,.
\end{eqnarray}
If we show that, for any~$1\leq s\leq t$, there is a~$\nu_s\in\mathbb N$ such that
\begin{eqnarray}\label{tesina}
B_s\subset\Omega\setminus\overline{\im V((F_{\chi_o})_j)}\,,
\end{eqnarray}
when~$j\geq\nu_s$, it will be enough to take~$k_1=\max\{\nu_1,\ldots,\nu_t\}$. So let us prove~(\ref{tesina}) by contradiction. 
\par
Suppose that, for some~$1\leq s\leq t$, there is an increasing sequence of integers~$j_q$ such that~$B_s$ intersects the closure of~$\im V(F_{\chi_o})_{j_q}$, (which, in particular, requires~$\im V(F_{\chi_o})_{j_q}\neq\varnothing$), for every integer~$q\in\mathbb N$. For any such~$q$, let~$\eta_q$ be a point in~$B_s\cap\overline{\im V(F_{\chi_o})_{j_q}}$, and, by compactness of~$B_s$, let~$q_\iota$ be an increasing sequence of integers such that the subsequence~$(\eta_{q_\iota})_\iota$ converges to some~$\eta\in B_s$. For each fixed~$\iota$, let also~$(\beta_{q_\iota,\ell})_\ell\subset\im V(F_{\chi_o})_{j_{q_\iota}}$ be a sequence that converges to~$\eta_{q_\iota}$, as~$\ell\to+\infty$, and~$(\alpha_{q_\iota,\ell})_\ell\subset\mathbb R^n$ a sequence such that
\begin{eqnarray}\label{equ1}
(F_{\chi_o})_{j_{q_\iota}}(\alpha_{q_\iota,\ell}+i\beta_{q_\iota,\ell})=0\,,
\end{eqnarray}
for every~$\iota,\ell\in\mathbb N$. If, for every fixed~$\iota$ and~$\ell$ in~$\mathbb N$,~$\tilde\chi_{q_\iota,\ell}$ denotes the continuous character of~$\mathbb R^n$ given by~$\exp(i\langle\alpha_{q_\iota,\ell},\cdot\rangle)$, equation~(\ref{equ1}) simply reads
\begin{eqnarray}\label{equ2}
((F_{\chi_o})_{j_{q_\iota}})_{\tilde\chi_{q_\iota,\ell}}(i\beta_{q_\iota,\ell})=0\,,
\end{eqnarray}
for every~$\iota,\ell\in\mathbb N$. Set~$G=F_{\chi_o}$ and observe that for any fixed~$\iota$, the sequence~$(\tilde\chi_{q_\iota,\ell})_\ell$ admit a convergent subsequence~$(\tilde\chi_{q_\iota,\ell_p})_p$ whose limit will be denoted~$\gamma_{q_\iota}$. This implies that
\begin{eqnarray}\label{equ3}
(G_{j_{q_\iota}})_{\tilde\chi_{q_\iota,\ell_p}}(i\beta_{q_\iota,\ell_p})=0\,,
\end{eqnarray}
for every~$\iota$ and~$p$, hence Proposition~\ref{facilis} implies
\begin{eqnarray}\label{equ4}
0
=
\lim_{p\to+\infty}
(G_{j_{q_\iota}})_{\tilde\chi_{q_\iota,\ell_p}}(i\beta_{q_\iota,\ell_p})
=
(G_{j_{q_\iota}})_{\gamma_{q_\iota}}(i\eta_{q_\iota})
\,,
\end{eqnarray}
for every~$\iota$. Now,~$(\gamma_{q_\iota})_\iota$ has a convergent subsequence~$(\gamma_{q_{\iota_p}})_p$, so let us call~$\gamma$ its limiting character. 
By Corollary~\ref{utile2},
\begin{eqnarray}\label{equ5}
\lim_{p\to+\infty}
(G_{j_{q_{\iota_p}}})_{\gamma_{q_{\iota_p}}}(i\eta_{q_{\iota_p}})
=
G_{\gamma}(i\eta)\,,
\end{eqnarray}
then~(\ref{equ4}) yields~$G_\gamma(i\eta)=0$, or equivalently~$F_{(\chi_o\,\cdot\,\gamma)}(i\eta)=0$, which is impossible since~$\eta\in B_s\subset\Omega\setminus{\mathcal F}_F$.
\par
\smallskip
\noindent
\textit{Step~4.} 
This Step is almost the same as the corresponding one in the proof of~Theorem~\ref{buono}. If~$j$ is bigger than~$k_1$ and~$y$ is a point of~$S\cap\overline{\im V(F_{\chi_o})_j}$, then
$$
\int_c\upsilon_y>0\qquad{\hbox{\rm if~$y\in\Supp C$,}}
\qquad
\hbox{\rm whereas}
\qquad
\int_c\upsilon_y=0\qquad{\hbox{\rm if~$y\not\in\Supp C$}}\,,
$$
so, by Stokes theorem, for any~$j\geq k_1$, the chain~$c$ represents a non negative homology class in the group~$H_m(S\setminus\overline{\im V((F_{\chi_o})_j)})$. Now, since~$\Omega$ is convex, the inclusion~$\Supp c\subset\Omega$ implies that~$\Omega$ contains the convex hull of~$\Supp c$. As~$\Supp C$ is included in the convex hull of~$\Supp c$, it follows that~$\Supp C\subset\Omega$, and, again by convexity,~$\Omega$ will contain the convex hull of~$\Supp C$. Moreover,~$\Omega$ is open, so it contains some~$\varepsilon$-neighbourhood~$\Omega^\prime$ of the convex hull of~$\Supp C$. The subset~$\Omega^\prime$ is non-empty, open, convex and relatively compact in~$\Omega$, hence (by Step~1) we may find a~$k_0\in\mathbb N$ such that~$\Omega^\prime\setminus\overline{\im V((F_{\chi_o})_j)}$ is~$m$-convex, for every~$j\geq k_0$. Observing that, for every~$j\geq k_1$, the chain~$c$ represents a non-negative class in the homology group~$H_m(S\cap\Omega^\prime\setminus\overline{\im V(F_{\chi_o})_j)})$, it follows (by $m$-convexity) that this class becomes trivial as soon as~$j\geq k_2=\max(k_0,k_1)$. Equivalently this means that
\begin{eqnarray*}\label{assurdo1}
\Supp C\subset S\cap\Omega^\prime\setminus\overline{\im V((F_{\chi_o})_j)}\,,
\end{eqnarray*}
when~$j\geq k_2$, so {\it a fortiori}
\begin{eqnarray*}
\Supp C\subset S\cap\Omega^\prime\setminus{\im V((F_{\chi_o})_j)}\,,
\end{eqnarray*}
for any~$j\geq k_2$.
\par
\smallskip
\noindent
\textit{Step~5.}
Let~$y\in S\cap\Omega^\prime\cap \im V(F_{\chi_o})$ be a point in the relative interior of~$\Supp C$. Then going on as in the Step~5 of the proof of~Theorem~\ref{buono}, we get, for any sufficiently big~$j$, a point~$y_j$ which should belong to the empty set, namely the intersection of~$S\cap\Omega^\prime\cap\im V((F_{\chi_o})_j)$ with~$S\cap\Omega^\prime\setminus\im V((F_{\chi_o})_j)$. The contradiction completes the proof. \hfill$\square$
\par
\begin{remark}
{\rm 
Corollary~\ref{buonis} is a substantial improvement of Theorem~5.1 in~\cite{sil06}, but in fact,~Theorem~\ref{buono} itself is already an improvement of that theorem. The next section will give a detailed account of the relationships between the hypothesis of~Theorem~\ref{buono} and those of Theorem~5.1 in~\cite{sil06}.}
\end{remark}
\par

\section{On regularity}\label{fin}
\subsection{Kazarnovski{\v\i} regularity}

In the paper~\cite{ka81}, Kazarnovski{\v\i} defined a notion of regularity for finite systems of exponential sums. In particular, this notion applies to the mappings belonging to the space~$\Ex(\mathbb C^n,\mathbb C^m)$. We recall the definition in our situation. 
\par
Let~$F=(f_1,\ldots,f_m)\in\Ex(\mathbb C^n,\mathbb C^m)$. Let us consider, for every~$1\leq\ell\leq m$, the {\sl Newton polytope}~$\Gamma_\ell$ of~$f_\ell$, that is the convex hull of the set
$$
\{-i\lambda\in i\mathbb R^n\mid\lambda\in\Sp f_\ell\}
$$
and let
$$
h_{\Gamma_\ell}(z)
=
\sup_{u\in\Gamma_\ell}\re\langle z,u\rangle
=
\sup_{\lambda\in\Sp f_\ell}\re\langle z,-i\lambda\rangle
=
-\sup_{\lambda\in\Sp f_\ell}\langle \im z,\lambda\rangle
$$
be its {\sl support function}. If~$\Gamma_F$ denotes the {\sl Minkowski sum} of the convex sets~$\Gamma_1,\ldots,\Gamma_m$, we recall that any face~$\Delta$ of~$\Gamma_F$ has a unique representation~$\Delta=\sum_{\ell=1}^m\Delta_\ell$, where, for every~$1\leq\ell\leq m$,~$\Delta_\ell$ is a face of~$\Gamma_\ell$. To any face~$\Delta$ of~$\Gamma_F$ one can associate a truncated mapping~$F^\Delta$ of~$F$, called the~$\Delta${\sl -trace} of~$F$. Its components~$f_1^\Delta,\ldots,f_m^\Delta$ are truncations of those of~$F$, namely, for every~$1\leq\ell\leq m$,
$$
f_\ell^\Delta(z)
=
\sum_{\lambda\in i\Delta_\ell\cap\Sp f_\ell}
a(\lambda,f_\ell)\,e^{i\langle z,\lambda\rangle}\,.
$$
\begin{remark}
{\rm
We observe that, for any character~$\chi\in\Ch\mathbb R^n$ and any face~$\Delta$ of~$\Gamma_F$, we have~$(F^\Delta)_\chi=(F_\chi)^\Delta$, hence we may simply write~$F_\chi^\Delta$ without warring about the order with which the truncation and the perturbation are performed.}
\end{remark}
\par
A mapping~$F\in\Ex(\mathbb C^n,\mathbb C^m)$ is said to be {\sl Ka\-zar\-nov\-ski{\v\i} regular}~(cf.~\cite{ka81}) if, for any face~$\Delta=\sum_{\ell=1}^m\Delta_\ell$ of~$\,\Gamma_F$ whose dimension~$\dim\Delta$ is strictly smaller than~$m$, there exists~$\varepsilon=\varepsilon(F,\Delta)>0$ such that
$$
K[F^\Delta](z)
=
\sum_{\ell=1}^m
e^{-h_{\Delta_\ell}(z)}
\Big\vert
f_\ell^\Delta(z)
\Big\vert
\geq\varepsilon\,,
$$
for every~$z\in\mathbb C^n$.
\par
\begin{remark}
{\rm 
If a mapping~$F\in \Ex(\mathbb C^n,\mathbb C^m)$, with~$m>n$, is Kazarnovski{\v\i} regular then it has empty zero set. In fact, if~$m>n$, the polytope~$\Gamma_F$ is a face of itself whose dimension is strictly smaller than~$m$, hence there exists an~$\varepsilon>0$ such that~$K[F^{\Gamma_F}](z)\geq\varepsilon$, for every~$z\in\mathbb C^n$; in particular, the function~$K[F^{\Gamma_F}]$ has empty zero set. Since~$F=F^{\Gamma_F}$ and~$F$ has the same zero set as~$K[F]$, the statement follows easily.}
\end{remark}
\par
\begin{lemma}{\rm (Kazarnovski{\v\i} \cite{ka81})}\label{kaza1}
Let~$F\in\Ex(\mathbb C^n,\mathbb C^m)$ be a Kazarnovski{\v\i} regular mapping. Then its zero set~$V(F)$ is a complete intersection.
\end{lemma}
\pf
An idea for the proof may be found in the original paper~\cite{ka81}. A complete proof is available in~\cite{sil05}.\hfill$\square$ 
\par
\begin{remark}\label{cav}
{\rm
As shown in the paper~\cite{ka81}, the zero set of a Kazarnovski{\v\i} regular mapping~$F\in\Ex(\mathbb C^n,\mathbb C^m)$,~$1\leq m\leq n$, is not empty if and only if the dimension of the affine subspace of~$\mathbb R^n$ generated by the union of~$\Sp f_1,\ldots,\Sp f_m$ is a least equal to~$m$ or, equivalently, if and only if the~$n$-dimensional mixed volume of~$i\Gamma_1,\ldots,i\Gamma_m$ with~$(n-m)$ full-dimensional balls of~$\mathbb R^n$ is strictly positive. We observe that both these latter properties are common to every mapping in the orbit~${\mathcal O}(F)$.}
\end{remark}
\par
\begin{lemma}\label{iokaza}
Let~$F=(f_1,\ldots,f_m)\in\Ex(\mathbb C^n,\mathbb C^m)$ be Ka\-zar\-nov\-ski{\v\i} regular. Then any mapping in the orbit~${\mathcal O}(F)$ is Kazarnovski{\v\i} regular too.
\end{lemma}
\pf
Let~$\chi\in\Ch\mathbb R^n$ and let~$\{\omega_1,\ldots,\omega_r\}$ a system of free generators of~$\Xi_F$. By~Poposition~\ref{krom}, there is a sequence~$(t_s)\subset\mathbb R^n$ such that
\begin{eqnarray}\label{kronek}
\lim_{s\to+\infty}
\langle t_s,\omega_j\rangle=\Arg\chi(\omega_j)\qquad\hbox{\rm mod.}2\pi\mathbb Z\,,
\end{eqnarray}
for any~$1\leq j\leq r$. Now if~$\Delta=\sum_{\ell=1}^m\Delta_\ell$ is face of~$\Gamma_F=\sum_{\ell=1}^m\Gamma_\ell$ whose dimension is strictly smaller than~$m$, the equalities~(\ref{kronek}) imply, (just like in the proof of~Lemma~\ref{figo}), that
$$
(F_\chi)^\Delta(z)
=
\lim_{s\to+\infty}
F^\Delta(z+t_s)\,,
$$
in the topology~$\tau(T_{\mathbb R^n})$, hence
\begin{eqnarray*}
K\big[(F^\Delta)_\chi\big](z)
&=&
\sum_{\ell=1}^m
e^{-h_{\Delta_\ell}(z)}
\Big\vert
(f^\Delta)_\chi(z)\Big\vert
\\
&=&
\sum_{\ell=1}^m
\exp\bigg[{\sup_{\lambda\in i\Delta_\ell\cap\Sp f_\ell}\langle\im z,\lambda\rangle}\bigg]
\Big\vert
(f^\Delta)_\chi(z)\Big\vert
\\
&=&
\lim_{s\to+\infty}
\sum_{\ell=1}^m
\exp\bigg[{\sup_{\lambda\in i\Delta_\ell\cap\Sp f_\ell}\langle\im (z+t_s),\lambda\rangle}\bigg]
\Big\vert
f^\Delta(z+t_s)\Big\vert
\\
&=&
\lim_{s\to+\infty}
\sum_{\ell=1}^m
e^{-h_{\Delta_\ell}(z+t_s)}
\Big\vert
f^\Delta(z+t_s)\Big\vert
\\
&=&
\lim_{s\to+\infty}
K\big[F^\Delta\big](z+t_s)\,,
\end{eqnarray*}
in the same topology. Since~$F$ is regular, there is an~$\varepsilon>0$ such that
$$
K[F^\Delta](z+t_s)\geq\varepsilon\,,
$$
for any~$z\in\mathbb C^n$ and any~$s\in\mathbb N$, so~$K[(F_\chi)^\Delta]\geq\varepsilon$ too. The lemma follows by the arbitrary choice of~$\Delta$ among the faces of~$\Gamma_F$ whose dimension is strictly smaller than~$m$.\hfill$\square$

\par
\begin{corollary}\label{ronkaza}
A Kazarnovski{\v\i} regular mapping is regular.
\end{corollary}
\pf  The corollary is a consequence of~Theorem~\ref{figo2} and~Remark~3, together with~Lemma~\ref{iokaza}.\hfill$\square$
\par
The class of regular mappings belonging to~$\Ex(\mathbb C^n,\mathbb C^m)$ is strictly larger than that of Kazarnovski{\v\i} regular ones. In order to give a better account of the situation, let us first recall the notion of proper regularity introduced by Rashkovskii~\cite{ras01}. 

A mapping~$F\in HAP(T_\Omega,\mathbb C^m)$ is said to be \textsl{properly regular} if, for any mapping~$G=(g_1,\ldots,g_m)\in\overline{{\mathcal B}(F)}$ the zero set~$V(G)$ is a proper intersection, (i.e. for any non empty subset~$I\subseteq\{1,\ldots,m\}$, the intersection of the zero sets~$V(g_j)$, as~$j$ runs in~$I$, is either empty or~$(\card I)$-codimensional).
\par
A properly regular mapping~$F\in HAP(T_\Omega,\mathbb C^m)$ is always regular, but the converse is generally false, unless~$m\leq 2$. 
\par
Let us now show that Corollary~\ref{ronkaza} may not be reversed. Consider the mapping~$F=(f_1,f_2)\in\Ex(\mathbb C^2,\mathbb C^2)$ given by
\begin{eqnarray*}
f_1(z)
&=&
e^{iz_1}+e^{iz_2}+e^{i(z_1+z_2)}+2\,,
\\
f_2(z)
&=&
3e^{iz_1}-e^{iz_2}-e^{i(z_1+z_2)}+4\,.
\end{eqnarray*}
Both components have the same Newton polytope and, if~$\Delta$ is the face of~$\Gamma_F$ given by the line segment~$[(-2i,-2i),(0,-2i)]$, the corresponding~$\Delta$-trace has the following components
\begin{eqnarray*}
f_1^\Delta(z_1,z_2)
&=&
e^{iz_2}+e^{i(z_1+z_2)}
\\
f_2^\Delta(z_1,z_2)
&=&
-e^{iz_2}-e^{i(z_1+z_2)}\,,
\end{eqnarray*}
hence~$V(F^\Delta)\neq\varnothing$ and consequently~$F$ cannot be Kazarnovski{\v\i} regular. On the other hand, if~$G=(g_1,g_2)\in\Ex(\mathbb C^2,\mathbb C^2)$ is the mapping with
\begin{equation}\label{mapG}
g_1(z_1,z_2)
=
2e^{iz_1}+3
\qquad
\hbox{\rm and}
\qquad
g_2(z_1,z_2)
=
e^{iz_2}-1\,,
\end{equation}
it is easily seen that~$G$ is (properly) regular. Moreover~$V(F_\chi)=V(G_\chi)$, for any~$\chi\in\Ch\mathbb R^2$, so~$F$ is an example of a (properly) regular mapping which is not Kazarnovski{\v\i} regular. 
\par
By the way, we remark that the mapping~$G$ defined by~(\ref{mapG}) is both properly regular and Kazarnovski{\v\i} regular.
\par
We still observe that the mapping~$H=(h_1,h_2,h_3)\in\Ex(\mathbb C^3,\mathbb C^3)$ given by
\begin{eqnarray*}
h_1(z)
&=&
e^{iz_1}+e^{iz_2}+e^{i(z_1+z_2)}+2\,,
\\
h_2(z)
&=&
(e^{iz_1}+e^{iz_2}+e^{i(z_1+z_2)}+2)(3e^{iz_1}-e^{iz_2}-e^{i(z_1+z_2)}+4)
\,,
\\
h_3(z)
&=&
e^{iz_3}-1\,,
\end{eqnarray*}
is not properly regular but it is Kazarnovski{\v\i} regular (and even with closed spectra, cf. Subsection~\ref{reggeo}). Finally, notice that, unlike proper regularity, Kazarnovsk{\v\i} regularity may fail to be transmitted to``submappings'', in fact~$H$ is Kazarnovski{\v\i} regular whereas~$(h_1,h_2)\in\Ex(\mathbb C^3,\mathbb C^2)$ is not such. This failure to transmit Kazarnovski{\v\i} regularity may even affect those mappings which are simultaneously properly regular and Kazarnovski{\v\i} regular: take~$h_1$ and~$h_3$ as above, but~$h_2(z)=3e^{iz_1}-e^{iz_2}-e^{i(z_1+z_2)}+4$.
\par
\subsection{Sufficient conditions for regularity}\label{reggeo}

A mapping~$F=(f_1,\ldots,f_m)\in \Ex(\mathbb C^n,\mathbb C^m)$ is said to have {\sl closed spectra}~(cf.~\cite{ka81}) if, for any face~$\Delta=\sum_{\ell=1}^m\Delta_\ell$ of the polytope~$\Gamma_F$ whose dimension~$\dim\Delta$ is strictly smaller than~$m$, there exists an~$\ell\in\{1,\ldots,m\}$ such that the face summand~$\Delta_\ell$ consists of a single point.
\par
This geometric condition  involves only the spectra of the components of the mapping and, once the number of spectra and the number of points in each spectrum are fixed, the condition is generically fulfilled.
\par
\begin{lemma}\label{kazacl}
A mapping~$F=(f_1,\ldots,f_m)\in \Ex(\mathbb C^n,\mathbb C^m)$ with closed spectra is Kazarnovski{\v\i} regular.
\end{lemma}
\pf
The proof is straightforward.\hfill$\square$
\par
\begin{remark}
{\rm
It is well-known that Lemma~\ref{kazacl} is not reversible. For example, the mapping~$F=(f_1,f_2)\in\Ex(\mathbb C^2,\mathbb C^2)$ with
$$
f_1(z)=e^{iz_1}+2e^{iz_2}+1\,,
\qquad
f_2(z)=e^{iz_1}+3e^{iz_2}-1\,,
$$
has not closed spectra but it is Kazarnovski{\v\i} regular.
}
\end{remark}

\medskip
Lemma~\ref{kazacl},~Corollary~\ref{ronkaza} and~Theorem~\ref{buono} imply that the following theorem already proved in~\cite{sil06}.
\begin{theorem}{\rm (Silipo \cite{sil06})}\label{vecchio}
If~$F=(f_1,\ldots,f_m)\in \Ex(\mathbb C^n,\mathbb C^{m+1})$,~$m+1\leq n$, has closed spectra, its amoeba complement~$\mathbb R^n\setminus{\mathcal F}_F$ is an~$m$-convex subset.
\end{theorem}
\par
In fact, the closedness of the set of the spectra of a mapping~$F$ belonging to~$\Ex(\mathbb C^n,\mathbb C^{m+1})$ is a property that is shared by all the mappings in the orbit~${\mathcal O}(F)$. Moreover, for sufficiently big~$j$ and uniformly with respect to~$\chi\in\Ch\mathbb R^n$, this property is also shared by all the mappings in the sequence~$(F_\chi)_j$ constructed as in the Step~1 of the proof of~Theorem~\ref{buono}, whence the possibility to prove~Theorem~\ref{vecchio} in the paper~\cite{sil06} without using~Lemma~\ref{figo} and~Lemma~\ref{iokaza}.
\par
\medskip
In order to complete the picture of the relationships between Ronkin's regularity and Kazarnovski{\v\i}'s one, we recall that, given a general holomorphic almost periodic mapping~$F$, Ronkin~\cite{ro90} found a sufficient condition on the spectra of the components of~$F$ implying its regularity. We introduce the necessary notation to state this condition but we do it in a way that is slightly different with respect to the original paper~\cite{ro90}. This will turn out to be useful in the proof of~Lemma~\ref{geom}.
\par
If~$F=(f_1,\ldots,f_m)\in HAP(\mathbb C^n,\mathbb C^m)$, with~$1\leq m\leq n$, as before we consider, for every~$1\leq\ell\leq m$, the convex hull~$\Gamma_\ell$ of~$-i(\Sp f_\ell)$, together with its {\sl support function}
$$
h_{\Gamma_\ell}(z)
=
\sup_{v\in\Gamma_\ell}\re\langle z,v\rangle\,,
$$
(which may now take infinite values if~$\Gamma_\ell$ is unbounded). Let also, for~$u\in i\mathbb R^n$,
$$
H_{\Gamma_\ell}(u)
=
\{z\in i\mathbb R^n\mid \re\langle z,u\rangle= h_{\Gamma_\ell}(u)\}\,,
$$
and, for every~$\lambda\in\Lambda_\ell= \partial\Gamma_\ell\cap \Sp f_\ell$, consider the (possibly empty) cone
$$
K_{\lambda,\Gamma_\ell}
=
\{z\in i\mathbb R^n\mid \{\lambda\}=\Gamma_\ell\cap H_{\Gamma_\ell}(z)\}\,.
$$
Finally define
$$
{\mathcal Z}_\ell
=
i\mathbb R^n\setminus \bigcup_{\lambda\in \Lambda_\ell} K_{\lambda,\Gamma_\ell}\,,
\qquad
{\mathcal Z}_F
=
\bigcap_{\ell=1}^m{\mathcal Z}_\ell\,,
$$
and for every subset~$X$ of~$\mathbb R^n$, let~$\dim_H X$ denote the Hausdorff dimension of~$X$.
\begin{theorem}{\rm (Ronkin \cite{ro90})}\label{georonkin}
Let~$F\in HAP(\mathbb C^n,\mathbb C^m)$, with~$1\leq m\leq n$. The mapping~$F$ is regular as soon as~$\dim_H{\mathcal Z}_F\leq(n-m)$.
\end{theorem}
\par
\begin{remark}
{\rm 
As shown by the example given in the Remark~4 of the current section, the condition stated in the Theorem~\ref{georonkin} is not necessary for the regularity of the mapping. In fact, in that example, one has a regular mapping with~$\dim_H{\mathcal Z}_F=1>0=(n-m)$.
}
\end{remark}
The following lemma shows that, for a mapping belonging to~$\Ex(\mathbb C^n,\mathbb C^m)$, Kazarnovski{\v\i}'s geometric condition for regularity is equivalent to Ronkin's one.
\begin{lemma}\label{geom}
Let~$F\in\Ex(\mathbb C^n,\mathbb C^m)$, with~$1\leq m\leq n$. Then~$F$ has closed spectra if and only if~$\dim_H{\mathcal Z}_F\leq(n-m)$.
\end{lemma}
\pf
Let~$1\leq\ell\leq m$ be fixed and, for any face~$\Delta_\ell$ of~$\Gamma_\ell$ let us write~$\Delta_\ell\preccurlyeq\Gamma_\ell$. For any such face, consider the dual cone
$$
K_{\Delta_\ell,\Gamma_\ell}
=
\{z\in i\mathbb R^n\mid \Delta_\ell=\Gamma_\ell\cap H_{\Gamma_\ell}(z)\}\,.
$$
Observe that~$\Lambda_\ell$ contains the set~$V_\ell$ of the vertices of~$\Gamma_\ell$ and that~$K_{\lambda,\Gamma_\ell}=\varnothing$, for every~$\lambda\in\Lambda_\ell\setminus V_\ell$, so
\begin{eqnarray*}
{\mathcal Z}_\ell
&=&
i\mathbb R^n
\setminus 
\bigcup_{v\in V_\ell}
K_{v,\Gamma_\ell}
\\
&=&
\bigg(
\bigcup_{v\in V_\ell}
{\overline K}_{v,\Gamma_\ell}
\bigg)
\setminus 
\bigg(
\bigcup_{v\in V_\ell}
K_{v,\Gamma_\ell}
\bigg)
\\
&=&
\bigcup_{v\in V_\ell}
\bigg(
{\overline K}_{v,\Gamma_\ell}
\setminus
\bigcup_{v\in V_\ell}
K_{v,\Gamma_\ell}
\bigg)
\\
&=&
\bigcup_{v\in V_\ell}
\bigcup_{\Delta_\ell\preccurlyeq\Gamma_\ell\atop\Delta_\ell\neq\{v\}}
K_{\Delta_\ell,\Gamma_\ell}
\\
&=&
\bigcup_{\Delta_\ell\preccurlyeq\Gamma_\ell\atop\dim\Delta_\ell>0}
K_{\Delta_\ell,\Gamma_\ell}\,.
\end{eqnarray*}
It follows that
\begin{equation}\label{haus}
{\mathcal Z}_F
=
\bigcap_{\ell=1}^m
\bigcup_{\Delta_\ell\preccurlyeq\Gamma_\ell\atop\dim\Delta_\ell>0}
K_{\Delta_\ell,\Gamma_\ell}
=
\bigcup
\bigcap_{\ell=1}^m
K_{\Delta_\ell,\Gamma_\ell}\,,
\end{equation}
where the last union is taken over the set of~$m$-tuples~$(\Delta_1,\ldots,\Delta_m)$ of faces with $\Delta_\ell\preccurlyeq\Gamma_\ell$ and~$\dim\Delta_\ell>0$, for any~$1\leq \ell\leq m$.
For each such~$m$-tuple, the intersection on the right-hand side of~(\ref{haus}) provides a cone that spans an affine subspace having a dimension equal to the Hausdorff dimension of the cone itself, moreover, the whole~${\mathcal Z}_F$ is a star-shaped set consisting of the union of such cones, so
$$
\dim_H{\mathcal Z}_F
=
\dim
\bigg(
\bigcup
\bigcap_{\ell=1}^m
K_{\Delta_\ell,\Gamma_\ell}
\bigg)
=
\max
\bigg(\dim
\bigg(
\bigcap_{\ell=1}^m
K_{\Delta_\ell,\Gamma_\ell}
\bigg)
\bigg)
\,,
$$
where the maximum is taken on the set of~$m$-tuples~$(\Delta_1,\ldots,\Delta_m)$ with $\Delta_\ell\preccurlyeq\Gamma_\ell$ and~$\dim\Delta_\ell>0$, for any~$1\leq \ell\leq m$. Observe that, for such~$m$-tuples of faces, the dimension of the intersection
$$
\bigcap_{\ell=1}^m
K_{\Delta_\ell,\Gamma_\ell},
$$
is maximal when~$\Delta=\sum_{\ell=1}^m\Delta_\ell$ is a face of~$\Gamma_F=\sum_{\ell=1}^m\Gamma_\ell$, in which case
$$
K_{\Delta,\Gamma_F}
=
\bigcap_{\ell=1}^m 
K_{\Delta_\ell,\Gamma_\ell}\,.
$$
This means that the Hausdorff dimension of~${\mathcal Z}_F$ equals the maximum of the dimensions of the cones which are dual to the faces of~$\Gamma_F$ without zero-dimensional face summands. Since, for any~$\Delta\preccurlyeq\Gamma_F$,~$\dim K_{\Delta,\Gamma_F}=n-\dim\Delta$, it follows that the Hausdorff dimension of~${\mathcal Z}_F$ is equal to~$n$ minus the minimum of the dimensions of the faces of~$\Gamma_F$ which do not admit zero-dimensional face summands.
\par
This computation of~$\dim_H{\mathcal Z}_F$ implies the lemma. In fact, if the Hausdorff dimension of~${\mathcal Z}_F$ is smaller or equal to~$(n-m)$, the dimension of every face of~$\Gamma$ without zero-dimensional face summands will be bigger or equal to~$m$, which is the same as saying that~$F$ has closed spectra. On the other hand, if~$F$ has closed spectra, the minimum of the dimensions of the faces of~$\Gamma_F$ which do not admit zero-dimensional face summands is bigger or equal to~$m$, so that the Hausdorff dimension of~${\mathcal Z}_F$ is smaller or equal to~$(n-m)$.\hfill$\square$

\section*{Acknowledgements}
\par\noindent
The authors are grateful to Sergeij Yu. Favorov, Alexander Yu. Rashkovskii and Boris Ya. Kazarnovski{\v\i} for several valuable discussions, as well as to the anonymous referee for constructive criticism.
\par
%\vfill\eject

\bigskip
\noindent
\vfill\eject


\begin{thebibliography}{9}

\bibitem{ap71}%1 
L.~Amerio, G.~Prouse, 1971, {\it Almost-Periodic functions and functional equations,} Van Nostrand company. 

\bibitem{fa01}%2
S.~Favorov, 2001 Holomorphic almost periodic functions in tube domains and their amoebas, \textit{Comput. Methods Funct. Theory} {\bf 1}, No. 2, 403-415.

\bibitem{gkz94}%3
I.M.~Gelfand, M.M.~Kapranov, A.V.~Zelevinsky, 1994, {\it Discriminants, Resultants and multidimensional Determinants,} Birkha{\" u}ser, Boston.

\bibitem{hen04}%4
A.~Henriques, 2004, An analogue of convexity for complements of amoebas of varieties of higher codimension, an answer to a question asked by B.~Sturmfels, \textit{Adv. in Geom.} \textbf{4}, no. 1, 61-73.

\bibitem{ro90}%5
L.~I Ronkin, 1990, Jessen's theorem for holomorphic almost periodic mappings, \textit{Ukrainsk. Math. Zh.} \textbf{42}, 1094-1107.

\bibitem{sil06}%6
J.~Silipo, 2006, Amibes de sommes d'exponentielles, arXiv:math. CV/0401150 v3. To appear in \textit{Canad. J. of Math.}

\bibitem{mi04}%7
G.~Mikhalkin, 2004, Amoebas of Algebraic Varieties and Tropical Geometry, arXiv:math. AG/0403015 v1.

\bibitem{ras01}%8
A.~Yu.~Rashkovskii, 2001, {\sl Zeros of holomorphic almost periodic mappings with independent components,\/} Complex Var. vol~44, 299-316. 

\bibitem{sil05}%9
J.~Silipo, 2005, Systèmes de sommes d'exponentielles à spectres réels et structure de leurs amibes, PhD thesis, Université Bordeaux I, France.

\bibitem{y79}%10
A. Yger, 1979, Fonctions définies dans le plan et moyennes en tout point de leurs valeurs aux sommets de deux carrés, \textit{C.R. Acad. Sci. Paris}, Sér. A, {\bf 288}, no. 10, A535-A538.

\bibitem{yg79}%11
A. Yger, 1979, Une généralisation d'un théorème de J. Delsarte, \textit{C.R. Acad. Sci. Paris}, Sér. A, {\bf 288}, no. 9, A497-A499.

\bibitem{ka81}%12
B.Ya.~Kazarnovski{\v\i}, 1981, {\sl On the zeros of exponential sums,\/} Soviet Math.~Dokl. Vol. 23, no. 2, 347-351.


\end{thebibliography}
\end{document}